\magnification=1100
\font\bigbf=cmbx10 scaled \magstep2
\font\medbf=cmbx10 scaled \magstep1
\hfuzz=10 pt
\def\n{\noindent}

\def\R{{I\!\!R}}
\def\L{{\bf L}}

\def\sgn{{\rm sign}}
\def\vp{\varphi}
\def\M{{\cal M}}
\def\F{{\cal F}}

\def\C{{\cal C}}
\def\vp{\varphi}

\def\c{\centerline}
\def\wto{\rightharpoonup}
\def\meas{\hbox{meas}}
\def\v{\vskip 1em}
\def\vs{\vskip 2em}
\def\vsk{\vskip 3em}
\def\ve{\varepsilon}
\null
\vs
\c{\bigbf Global solutions of the Hunter-Saxton equation}
\v
\c{\it Alberto Bressan}
\v
\c{Deptartment of Mathematics, Pennsylvania State University,
University Park 16802, U.S.A.}
\c{e-mail: bressan@math.psu.edu}
\v
\c{and}
\v
\c{\it Adrian Constantin}
\v
\c{Trinity College, Department of Mathematics, Dublin 2, Ireland}
\c{e-mail: adrian@maths.tcd.ie}

\vsk
{\bf Abstract.} {\it
We construct a continuous semigroup of
weak, dissipative solutions to a
nonlinear partial differential equations modeling nematic liquid crystals.
A new distance functional, determined by
a problem of optimal transportation,
yields sharp estimates on the continuity of solutions with respect to the
initial data.
}

\vsk
\n{\medbf 1 - Introduction}
\v
In this paper we investigate the Cauchy problem
$$u_t+\left(u^2\over 2\right)_x
={1\over 4} \left( \int_{-\infty}^x-\int_x^\infty\right)
u_x^2\,dx\,,\qquad\qquad u(0,x)=\bar u(x)\,.\eqno(1.1)$$
Formally differentiating the above equation with respect to the spatial variable $x$, we obtain
$$(u_t+uu_x)_x={1 \over 2}\,u_x^2,\eqno(1.2)$$
whereas yet another differentiation leads to
$$u_{txx}+2u_xu_{xx}+uu_{xxx}=0.\eqno(1.3)$$
Either of the forms (1.2) and (1.3) of the equation in (1.1) is known as the
Hunter-Saxton equation.
In this paper we analyse various concepts of solutions
for the above equations, and construct a semigroup of globally defined
solutions.  Moreover, we introduce a new distance functional,
related to a problem of optimal transportation,
which monitors the continuous dependence of solutions on the initial data.

\bigskip

\noindent
{\bf Physical significance}\bigskip

The Hunter-Saxton equation  describes the propagation of waves in a
massive director field of a
nematic liquid cristal [HS], with the orientation of the molecules
described by the field of unit
vectors ${\bf n}(t,x)=\big(\cos\,u(t,x),\, \sin \,u(t,x)\big)$, $x$ being the
space variable
in a reference frame moving with the linearized wave velocity, and $t$ being
a slow time variable.
The liquid crystal state is a distinct phase of matter observed between the
solid and liquid states.
More specifically, liquids are isotropic (that is, with no directional order)
and without a positional order of their molecules, whereas the molecules
in solids are constrained to point only in certain directions and to be only in
certain positions with respect to each other. The liquid crystal phase exists
between the solid and the liquid phase - the molecules in a liquid crystal do
not exhibit any positional order, but they do possess a certain degree of
orientational order. Not all substances can have a liquid
crystal phase e.g. water molecules melt directly from solid crystalline ice to liquid water.
Liquid crystals are fluids made up of long rigid molecules, with
an average orientation that specifies the local direction of the medium.
Their orientation is described
macroscopically by a field of unit vectors ${\bf n}(t,{\bf x}))$. There are
many types of liquid crystals, depending upon the amount of order in the material.
Nematic liquid crystals are invariant
under the transformation ${\bf n} \mapsto -{\bf n}$, in  which case ${\bf n}$
is called a director field, so
that the rodlike molecules have no positional order but tend to point in the
same direction (along the director). The director field does not remain
the same but generally fluctuates.
Obtaining the equation governing
the director field represents the crucial point for the modeling of nematic
liquid crystals since it is
advantageous to study the dynamics of director field instead of studying the
dynamics of all the molecules. The fluctuations
of the director field are mainly due to the thermodynamical force caused by
elastic deformations in the form of twisting, bending, and splay (the latter
being a fan-shaped spreading out from the original direction, bending being
a change of direction, while twisting corresponds to a rotation of the
direction in planes orthogonal to the axis of rotation).
Consider director fields that lie on a circle and depend on a single spatial
variable $x$ so that twisting is not allowed. To describe the dynamics of
the director field independently of the coupling
with the fluid flow, let $u(t,x)$ be the perturbation about a constant value.
The asymptotic equation for weakly
nonlinear unidirectional waves is precisely equation (1.2), obtained as the
Euler-Lagrange equation
of the variational principle
$$\delta \int_{t_1}^{t_2} \int_{\R} (u_tu_x+uu_x^2)\,dxdt=0,$$
for the internal stored energy of deformation of the director field if
dissipative effects are neglected
(corresponding to the case when inertia effects dominate
viscosity) - see [HS] for the details of
the derivation. Unlike other studies, in the Hunter-Saxton model the
kinetic energy of the director
field is not neglected. In the asymptotic regime in which (1.2)
is derived, the nondimensionalized
kinetic energy density is $u_x^2$ so that the condition
$$\int_{\R} u_x^2(t,x)\,dx < \infty \eqno(1.4)$$
has to hold at any fixed time $t$ for a physically meaningful solution to
the Hunter-Saxton equation.

Equation (1.1) is also relevant in other physical situations, e.g. it  is a
high-frequency limit of the
Camassa-Holm equation [DP], a nonlinear shallow water equation
[CH, J] modeling solitons [CH] as well as breaking waves [CE].

\bigskip
\noindent
{\bf Geometric interpretation}
\bigskip

An interesting aspect of the Hunter-Saxton equation (see [KM])
is the fact that,
for spatially periodic functions, it describes geodesic flow on the homogeneous space
$\hbox{Diff}({\cal S})/\hbox{Rot}({\cal S})$ of the
infinite-dimensional Lie group $\hbox{Diff}({\cal S})$
of smooth orientation-preserving diffeomorphisms of the
unit circle ${\cal S}$ modulo the rotations
$\hbox{Rot}({\cal S})$, with respect to the right-invariant
homogeneous metric
$\langle f,\,g \rangle=\displaystyle\int_{\cal S}f_xg_x\,dx$.
The geometric interpretation of the
Hunter-Saxton equation establishes a natural
connection with the Camassa-Holm equation, which describes geodesic
flow on $\hbox{Diff}({\cal S})$
with respect to the right-invariant metric $\langle f,\,g \rangle=
\displaystyle\int_{\cal S}(fg+f_xg_x)\,dx$, see
[K, CK]. A similar geometric interpretation of (1.1) on the
diffeomorphism group of the line holds also
for smooth initial data $\bar u$ in certain weighted function spaces
but the involved technicalities are
more intricated (see [C] for the case of the
Camassa-Holm equation).

\bigskip
\noindent
{\bf Integrable structure}
\bigskip

The Hunter-Saxton equation has a an integrable structure.
The equation has a reduction (see [BSS, HZ1])
to a finite dimensional completely integrable Hamiltonian system
whose phase space consists of piecewise
linear solutions of the form
$$u(t,x)=\sum_{i=1}^n \alpha_i(t)\,|x-x_i(t)|,\eqno (1.5)$$
with the constraint
$$\sum_{i=1}^n \alpha_i(t)=0,\eqno(1.6)$$
the Hamiltonian being
$$H(x,\alpha)={1 \over 2}\,\sum_{i,j=1}^n \alpha_i\alpha_j\,|x_i-x_j|.$$
Due to their lack of regularity, functions of the form (1.5) are not
classical solutions of (1.2).
Below we will discuss in what sense they are weak solutions of the
Hunter-Saxton equation.
Let us point out that the constraint (1.6) is the necessary and sufficient
condition to ensure
that the distributional derivative $x \mapsto u_x(t,x)$ of a function
of the form (1.5) belongs
to the space $L^2(\R)$. Thus (1.4) holds.

In the family of smooth functions $u:\R\mapsto\R$  all of whose derivatives
$\partial_x^n u$ decay rapidly
as $x\to\pm\infty$, the Hunter-Saxton equation is
bi-Hamiltonian [HZ1]. If $D^{-1}$ is the skew-adjoint anti-derivative
operator given by
$$(D^{-1}f)(x)={1\over 2} \left( \int_{-\infty}^x-\int_x^\infty\right)
f(x)\,dx\,,\qquad\qquad f \in {\cal D}(\R),$$
the first Hamiltonian form for the Hunter-Saxton equation is
$$u_t=J_1\,{\delta{\cal H}_1 \over {\delta u}},
\qquad J_1=u_xD^{-2}-D^{-2}u_x,\qquad {\cal H}_1
={1 \over 2}\,\int_{\R} u_x^2\,dx,$$
whereas the second, compatible, Hamiltonian structure is
$$u_t=J_2\,{\delta{\cal H}_2 \over {\delta u}},\qquad
J_2=D^{-1},\qquad {\cal H}_2={1 \over 2}\,\int_{\R} uu_x^2\,dx,$$
Moreover, the Hunter-Saxton equation is formally integrable e.g. it has an
associated Lax pair (see [BSS]). However, the complete integrability of the
equation has been established only in the previously mentioned case when it
reduces to a finite dimensional
dynamical system.

\bigskip
\noindent
{\bf The notion of solution}
\bigskip

Physically relevant solutions of the Hunter-Saxton equation need to be of
finite kinetic energy so
that (1.4) must hold. This leads naturally to functions $u(t,x)$ with
distributional derivative
$u_x(t,\cdot)$ square integrable at every instant $t$. Note that the
integrability assumption $u_x(t,\cdot) \in L^2(\R)$ already imposes a
certain degree of regularity on the function $u$.
This suggests that it might be possible to incorporate a reasonably
high degree of regularity
in the concept of weak solution to the  Hunter-Saxton equation. Let us
first consider the concept
of weak solutions introduced by Hunter and Zheng [HZ2].\bigskip

 {\bf Definition 1.1}  {\it A function $u(t,x)$ defined on $[0,T]
\times \R$ is a solution
 of the equation (1.2) if

(i) $u \in C([0,T] \times \R;\R)$ and  $u(0,x)=\bar u(x)$ pointwise on $\R$;

(ii) For each $t \in [0,T]$, the map $x \mapsto u(t,x)$ is absolutely
continuous with $u_x(t,\cdot) \in L^2(\R)$. Moreover, the map
$t \mapsto u_x(t,\cdot)$ belongs to the space $L^\infty([0,T]; L^2(\R))$
and is locally Lipschitz continuous on $[0,T]$ with values in
$H^{-1}_{loc}(\R)$;

(iii) Equation (1.2) holds in the sense of distributions.}\bigskip

\noindent
Here and below, by a mapping $f$ that is locally Lipschitz or
locally bounded on $[0,T]$ with values in
$H^{-1}_{loc}(\R)$ we understand the following: for every $n \ge 1$
there is a constant $K_n \ge 0$
such that
$$\sup_{\{\psi \in {\cal D}(-n,n):\ \| \psi \|_{H^1(\R)} \le 1\}}\Big|
\langle f(t)-f(s), \psi \rangle\Big| \le K_n\,|t-s|,\qquad t,s \in [0,T],$$
respectively
$$\sup_{\{\psi \in {\cal D}(-n,n):\ \| \psi \|_{H^1(\R)} \le 1\}} \Big|
\langle f(t), \psi \rangle\Big| \le K_n,\qquad t \in [0,T].$$
Here ${\cal D}(a,b)$ is the family of smooth functions $f: \R \to \R$
with compact support within $(a,b) \subset \R$.

To a function $u: [0,T] \times \R \to \R$ with the above properties
associate the function
$F: [0,T] \times \R \to \R$ defined by
$$F(t,x)={1\over 4} \left( \int_{-\infty}^x-\int_x^\infty\right)u_x^2\,dx.
\eqno(1.7)$$
Then $F \in L^\infty_{loc}([0,T] \times \R;\R) \subset L^2_{loc}([0,T]
\times \R;\R)$.
Moreover, $F_x={1 \over 2}\,u_x^2$ so that equation (1.2) becomes
$$(u_t+uu_x-F)_x=0 \eqno(1.8)$$
in the sense of distributions. Note that $uu_x \in L^2_{loc}([0,T]
\times \R;\R)$. From
(1.8) we infer the existence of a distribution $h(t)$ so that
$u_t+uu_x-F=h(t) \otimes 1(x)$,
where $1(x)$ stands for the constant function with value $1$ on
$\R$. If $H(t)$ is a primitive
of the distribution $h(t)$, we deduce that the distribution
$U=u - H(t) \otimes 1(x)$ satisfies
$U_t=u_t - h(t) \otimes 1(x)=F-uu_x   \in L^2_{loc}([0,T] \times \R;\R)$
and
$U_x=u_x \in L^2_{loc}([0,T] \times \R;\R)$.
Therefore $U \in H^1_{loc}([0,T] \times \R)$.
Moreover, since $U_t=F-uu_x \in L^\infty_{loc}([0,T];H^{-1}_{loc}(\R))$
ensures that $U$ is
locally Lipschitz as a function from $[0,T]$ to $H^{-1}_{loc}(\R)$ and
so is also $u$, we deduce
that $h(t) \otimes 1(x) =u-U$ shares this property too. But then $h: [0,T]
\to \R$ has to be Lipschitz
continuous. We infer that $u=U+H(t) \otimes 1(x)$ belongs to the space
$H^1_{loc}([0,T] \times \R)$.
Since the requirement (iii) in Definition 1.1 ensures that the identity
$$\int_0^T\int_{\R} \Bigl( \phi_{xt}u + {1 \over 2}\,
\phi_{xx}u^2-{1 \over 2}\,\phi u_x^2\Bigr)\,dxdt=0$$
holds for every smooth function $\phi:(0,T) \times \R \to \R$
with compact support  in $(0,T) \times \R$,
we see that the notion of weak solution in the sense of
Definition 1.1 is stronger than the concept
of weak solution introduced by Hunter and Saxton [HS].
Another useful conclusion that can be
drawn from the previous considerations is that for a function $u$
with regularity properties
specified in (i)-(ii) of Definition 1.1, the requirement (iii) from
Definition 1.1 is equivalent to asking
that the equation
$$u_t+uu_x=F+h(t) \circ 1(x)\eqno(1.9)$$
holds in distribution sense for some Lipschitz continuous function
$h: [0,T] \to \R$. Any such function
$h$ is admissible. Among all these possibilities the most natural
one corresponds to the special
choice $h \equiv 0$. This leads us to the form (1.1) of the
Hunter-Saxton equation.

In the following, we say that a map $t\mapsto u(t,\cdot)$ from $[0,T]$
into $\L^p_{loc}(\R)$ is {\it absolutely continuous} if, for every
bounded interval $[a,b]$, the restriction of $u$ to $[a,b]$ is
absolutely continuous as a map with values in $\L^p\big([a,b]\big)$.
We can thus adopt
the following notion of a weak solution.\bigskip

{\bf Definition 1.2} {\it A function $u(t,x)$ defined on
$[0,T] \times \R$ is a solution
 of the equation (1.2) if

(i) $u \in C([0,T] \times \R;\R)$ and  $u(0,x)=\bar u(x)$
pointwise on $\R$;

(ii) For each $t \in [0,T]$, the map $x \mapsto u(t,x)$ is
absolutely continuous with $u_x(t,\cdot) \in L^2(\R)$.
Moreover, the map $t \mapsto u_x(t,\cdot)$ belongs to
the space $L^\infty([0,T]; L^2(\R))$;

(iii) The map $t \mapsto u(t,\cdot) \in L^2_{loc}(\R)$ is
absolutely continuous and satisfies
the equation (1.1) for a.e. $t \in [0,T]$.}\bigskip

\noindent
The concept of solution introduced in Definition 1.2 is stronger than
that corresponding
to Definition 1.1. Indeed, for a function $u$ satisfying all the
requirements of Definition 1.2 we
infer by (1.1) that $u_{tx} \in L^\infty_{loc}([0,T];H^{-1}_{loc}(\R))$
since $u_t=-uu_x+F$ and
$uu_x,\,F \in L^\infty_{loc}([0,T];L^2_{loc}(\R))$.
This yields that the map $t \mapsto u_x(t,\cdot)$ is locally
Lipschitz continuous on $[0,T]$ with values in $H^{-1}_{loc}(\R)$.
We thus recover
the apparently missing part from the requirement (ii) in Definition 1.1.

We remark that, even with this stronger definition, solution
are far from unique.  For example, consider the initial data
$${\bar u}(x)=0.\eqno(1.10)$$
There are now two ways to prolong the solution for times $t>0$.
On one hand, we can define
$$u(t,x)=0\qquad\qquad x\in\R\,, ~~t\geq 0\,.\eqno(1.11)$$
On the other hand, the function
$$u(t,x)\doteq\cases{ -2t\qquad &if\quad $x\leq -t^2$\cr
&\cr
{2x\over t}\qquad &if\quad $|x|< t^2$\cr
&\cr
2t\qquad &if\quad $x\geq t^2$\cr}\qquad\hbox{for}~~t\geq 0
\eqno(1.12)$$
provides yet another solution.
To distinguish between these two solutions, we need to consider
the evolution equation satisfied by the ``energy density'' $u_x^2$, namely
$$(u_x^2)_t+(uu_x^2)_x=0\,.\eqno(1.13)$$
For smooth solution, the conservation law (1.13) is satisfied pointwise.
Notice that the solution defined by (1.10), (1.12)
satisfies the additional conservation law
(1.13) in distributional sense, i.e.
$$\int\int_{\R_+ \times \R} \big\{u_x^2\vp_t + uu_x^2\,\vp_x\big\}\,dxdt=0
\eqno(1.14)$$
for every test function $\vp\in\C^1_c(\R_+ \times \R)$ whose compact support is contained
in the half plane where $t>0$.  On the contrary, the solution
defined by (1.10)-(1.11) dissipates energy.  More precisely, for every
$t_2\geq t_1\geq 0$ we have
$$\int_\R u_x^2(t_2,x)\,\vp(t_2,x)\,dx-
\int_\R u_x^2(t_2,x)\,\vp(t_2,x)\,dx
\leq \int_{t_1}^{t_2} \int_\R \big\{ u_x^2\,\vp_t+uu_x^2 \vp_x\big\}\,dxdt\,,
\eqno(1.15)$$
for every test function $\vp\in\C^1_c(\R_+ \times \R)$.
In the sequel, we say that a solution is {\it dissipative} if the
inequality (1.15) holds for every $t_2>t_1>0$, $\vp\in\C^1_c(\R_+ \times \R)$.
Notice that the solution (1.10), (1.12) does not satisfy (1.15)
when $t_1=2$, $t_2>0$.

At this point in the discussion it is worthwhile to point out that
the most important
feature in the definition of weak solutions is the requirement (1.4).
A continuous
function $u: [0,T] \times \R \to \R$ with square integrable distributional
derivative $u_x(t,\cdot)$
belonging to the space $L^\infty([0,T];L^2(\R))$ is not necessarily bounded,
nor does it
have a pre-determined asymptotic behaviour at infinity, as one can see
from the
example
$$u(t,x)=\cases{|x|^{1 \over 5}\,\sin(|x|^{1 \over 5})\qquad &if\qquad
$t \ge 0,\ |x| \ge 1\,$,\cr
|x|^{2 \over 3}\,\sin(|x|^{2 \over 3})\qquad &if\qquad
$t \ge 0,\ |x| \le 1\,$.\cr}$$
Nevertheless, the possibility that some additional
structural information about the behaviour
of such functions at infinity might be inferred from some
invariance properties of
the Hunter-Saxton equation should be ruled out.
To do this, consider solutions of the type
(1.5) with the constraint (1.6). This type of
solutions enter into the framework of
Definition 1.2 and for any $N(t) > \max\,\{|x_1(t)|,...,|x_n(t)|\}$
we have
$$u_t(t,x)=F(x)\qquad\hbox{a.e. on}\qquad |x| \ge N(t),$$
so that for all $j \ge N(t)$,
$$u_t(t,j)-u_t(-j)=F(j)-F(-j)={1 \over 2} \, \int_{-j}^j u_x^2(t,x)\,dx
={1 \over 2} \, \int_{\R} u_x^2(t,x)\,dx.\eqno(1.16)$$
But the quantity $I=\displaystyle{1 \over 2}\,\int_{\R} u_x^2(t,x)\,dx$
is an
invariant (time-independent) cf. [HZ1, BSS]. Moreover, the special
form of the solutions
guarantees that at every fixed $t \ge 0$,
$$u_\infty(t) \doteq \lim_{x \to \infty}\,u(t,x)
=\sum_{i=1}^n \alpha_i(t)\,x_i(t)=-\lim_{x \to -\infty}\,u(t,x).$$
and $u_\infty(t)=u(t,j)=-u(t,-j)$ for
all $j \ge N(t)$. Thus (1.16) yields
$$u_\infty(t)=u_\infty(0)+{1 \over 4}\,It,\qquad t \ge 0.
\eqno(1.17)$$
Unless $I =0$, in which case $u$ is constant, we see from (1.17)
that the asymptotic behaviour
of the solutions changes with time. On the basis of this set
of examples we conclude that
the asymptotic behaviour of the solutions at infinity should
not be prescribed {\it a priori}. However, the previous set of
examples indicates that a possible
restriction would be to require $u \in L^\infty([0,T] \times \R)$
if $\bar u \in L^\infty(\R)$.
In this case the space of functions introduced in Definition 1.2
(that is,
bounded functions with all the properties specified in
Definition 1.2 except the
condition that $u$ satisfies equation (1.1) in $L^2[-n,n]$
for every $n \ge 1$) is
a Banach space when endowed with the norm
$$\| u \|_T=\sup_{(t,x) \in [0,T] \times \R}\,\{|u(t,x)|\}
+ \hbox{ess-}\!\sup_{\!\!\!\!\!\!\!\!\!\! t \in [0,T]}
\,\int_{\R} u_x^2(t,x)\,dx.\eqno(1.18)$$
It is also worth noticing that a function entering the framework of
Definition 1.2
has further reqularity properties that are not explicitely stated.
For example, we have the H\"older continuity property
$$|u(t,x)-u(t,y)| \le K(t)\,\sqrt{|x-y|},\qquad t \ge 0,\ x,y \in \R,$$
with $K(t)=\| u_x(t,\cdot) \|_{L^2(\R)}$, since
$$|u(t,x)-u(t,y)|^2 =\Bigl| \int_x^y u_x(t,\zeta)\,d\zeta
\Bigr|^2 \le |x-y|\,\cdot\, \Bigl|
\int_x^y u_x^2(t,\zeta)\,d\zeta\Bigr| \le |x-y|\,
\int_{\R} u_x^2(t,\zeta)\,d\zeta.$$

\vsk
\n{\medbf 2 - Global existence of dissipative solutions}
\v

For twice continuously differentiable solutions,
the derivative $v\doteq u_x$ of the
solution $u$ to (1.1) satisfies the equations
$$v_t+uv_x~=~-{v^2\over 2}\,,\eqno(2.1)$$
$$(v^2)_t+(uv^2)_x=0\,.\eqno(2.2)$$

Define the characteristic $t\mapsto \xi(t,y)$ as the solution to the O.D.E.
$${\partial\over\partial t}\,\xi(t,y)=u\big(t,\,\xi(t,y)\big)\,,
\qquad\qquad\xi(0,y)=y\,.\eqno(2.3)$$
From (1.2) it follows that the evolution of the gradient $u_x$
along each characteristic is
described by
$${d\over dt}\, u_x\big(t,\xi(t,y)\big)=-{1\over 2}
u^2_x\big(t,\xi(t,y)\big)\,.
\eqno(2.4)$$
Observe that the solution of the O.D.E.
$$\dot z=-z^2/2\,,\qquad\qquad z(0)=z_0$$
is given by
$$z(t)={2z_0\over 2+tz_0} \eqno(2.5)$$
If $z_0 \ge 0$ this solution is defined for all $t \ge 0$,
whereas if $z_0<0$, this solution
approaches $-\infty$ at the blow-up time
$$T(z_0)=-2/z_0\eqno(2.6)$$
Note that if $\bar u(x) \not \equiv 0$ then there is some $x_0 \in \R$
with $\bar u(x_0)<0$
so that the characteristic cruve $t \mapsto \xi(t,\bar u(x_0))$ will
blowup in finite time.
Nevertheless, if $\liminf_{x \in \R} \{\bar u_x(x)\}>-\infty$,
then $T_0>0$, where
$$T_0=\inf_{\{x \in \R:\ \bar u_x(x) <0\}}
\Big\{ {-2 \over \bar u_x(x)}\Big\} \ge 0,\eqno(2.7)$$
and on the time interval $[0,T_0)$ the method of
characteristics can be used to construct
the unique solution of (1.1). Let us describe the
construction in detail. From (2.3) we get
$$ {\partial\over\partial t}\,\xi_x=u_x(t,\xi)\,\cdot\,\xi_x=
{2\,\bar u_x \over 2+t\,\bar u_x}\,\cdot\,\xi_x
\eqno (2.8)$$
since
$$u_x(t,\xi(t,y))={2\,\bar u_x(y) \over 2+t\,\bar u_x(y)}\eqno(2.9)$$
in view of (2.4) and the solution formula (2.5). The unique solution
of the linear O.D.E. (2.7)
with initial data $\xi_x(0,y)=1$ is given by
$$\xi_x(t,y)=\Bigl(1+{t \over 2}\,\bar u_x(y) \Bigr)^2. \eqno(2.10)$$
Since $1+{t \over 2}\,\bar u_x(y)>0$ for $t \in [0,T_0)$, relation (2.10)
shows that for each
$t \in [0,T_0)$ the map $y \mapsto \xi(t,y)$ is an absolutely continuous
 increasing
diffeomorphism of the line. Define the absolutely continuous function
$\varphi$ by
$$\varphi(y)={1 \over 4}\,\int_{\R} \hbox{sign}(y-x)\,
\bar u_x^2(x)\,dx\eqno(2.11)$$
so that $\varphi_x(y)={1 \over 2}\,\bar u_x^2(y)$. Note that by (2.10),
$$\xi_{tx}=\bar u_x+t\,{\bar u_x^2 \over 2}.\eqno(2.12)$$
Since $\xi_t(0,y)=0$ as $\xi(0,y)=y$, integration of (2.12)
with respect to the spatial variable
$x$ yields
$$\xi_t(t,y)=\bar u(y)+{t \over 4}\, \int_{\R} \hbox{sign}(y-x)\,
\bar u_x^2(x)\,dx \eqno(2.13)$$
and thus
$$\xi(t,y)=y+\int_0^t \xi_t(s,y)\,ds=y+t\,\bar u(y)+
{t^2 \over 4}\,\int_{\R} \hbox{sign}(y-x)\,\bar u_x^2(x)\,dx.\eqno(2.14)$$
The value of the solution $u$ along the characteristic curve
$t \mapsto \xi(t,y)$ is
$$u(t,\xi(t,y))=\bar u(y) + {t \over 2}\, \int_{\R}
\hbox{sign}(y-x)\, \bar u_x^2(x)\,dx. \eqno(2.15)$$
This relation is obtained by combining (2.14) with (2.3).
The increasing diffeomorphism of the
line $y \mapsto \xi(t,y)$ given by (2.14) and formula (2.15)
yield the unique solution of
the Hunter-Saxton equation on the time interval $[0,T_0)$.
The above approach works as
long as $2+t\,\bar u_x(x)>0$ but breaks down at $T=T_0$
with $T_0$ given by (2.7). The
reason for the breakdown is that
$$\liminf_{t \uparrow T_0} \,\{ \inf_{x \in \R}\,u_x(t,x)\}=-\infty
\eqno(2.16)$$
in view of (2.9) and the definition (2.7) of $T_0$.
Note that at $t=T_0$ we have might have $\xi_x(t,x)=0$
for all $x \in (a,b) \subset \R$ so that the map $y \mapsto \xi(t,y)$ is not
any more an increasing diffeomorphism of the line.
Nevertheless, the previous considerations suggest the following
 approach in the general case when $\bar u_x \in L^2(\R)$, covering
situations when possibly
 $T_0=0$ as it is the case for e.g. $\bar u (x)=x^{2 \over 3}
(1-x)^{2 \over 3}\,\chi_{[0,1]}$. Here
 $\chi_A$ stands for the characteristic function of the set $A$,
defined by $\chi_A(x)=1$ if
 $x \in A$ and $\chi_A(x)=0$ if $x \not \in A$.

 Let $\bar u \in C(\R)$ be such that its distributional derivative
$\bar u_x$ is square integrable.
 Define $\varphi: \R_+ \times \R \to \R$ by
 $$\varphi_y(t,y)= {1 \over 2}\, \bar u_x^2(y)\,
\chi_{[\bar u_x > - 2/t]}(y) \eqno(2.17)$$
 so that
 $$\varphi(t,y)={1 \over 4}\, \int_{[\bar u_x > - 2/t]}
\hbox{sign}(y-x)\,\bar u_x^2(x)\,dx, \qquad
 t>0,\eqno(2.18)$$
 with the understanding that
 $$\varphi(0,y)={1 \over 4}\int_{\R}\hbox{sign}(y-x)\, \bar u_x^2(x)\,dx.$$
In other words, if $u(t,\xi(t,x_0)$ blows up before $t_0>0$, then
the point $x_0$ is not included in the domain of the integral defining
$\varphi(t_0,\cdot)$, because
 $$T\big(\bar u_x(x)\big)>t\qquad\hbox{if and only if}\qquad
\bar u_x(x)>-2/t\,,$$
according to (2.4) and (2.6). Observe that (2.18) and Young's
inequality yield
$$\| \varphi(t,\cdot) \|_{L^\infty(\R)} \le {1 \over 4}\,
\int_{\R} \bar u_x^2(x)\,dx,\qquad
t \ge 0.\eqno(2.19)$$

In the $(t,x)$-plane, the characteristic curve starting at $y$ is obtained as
$$\xi(t,y)= y+t\bar u(y)+\int_0^t (t-s)\,\varphi(s,y)\,ds\,.\eqno(2.20)$$
The value of the solution $u$ along this curve is
$$u\big(t,\,\xi(t,y)\big)=\bar u(y)+\int_0^t\varphi(s,y)\,ds\,.\eqno(2.21)$$
Observe that for all $t \ge 0$ and $y \in \R$,
$$\xi_t(t,y)= u(t,\xi(t,y))\eqno(2.22)$$
in view of (2.20)-(2.21).

\v
{\bf Theorem 1.} {\it Given any absolutely
continuous function $\bar u: \R \to \R$ with
derivative $\bar u_x \in L^2(\R)$, the formulas (2.18)-(2.20)
provide a dissipative
solution to (1.1), defined for all
times $t \ge 0$.}\bigskip

{\it Proof.} We proceed in several steps. First of all,
for any fixed $t \ge 0$, the map $y \mapsto \xi(t,y)$ is absolutely
continuous since
$\varphi_y(t,\cdot) \in L^2(\R)$. We claim that for any fixed $t \ge 0$
the map $y \mapsto \xi(t,y)$
is nondecreasing on $\R$ with $\lim_{y \to \pm\infty}\xi(t,y)=\pm\infty$.

Indeed, if $\bar u_x(y)>-\,{2 \over t}$ then $\bar u_x(y) >-\,{2 \over s}$
for all $s \in [0,t]$ so that
$\varphi_y(s,y)={1 \over 2}\,\bar u_x^2(y)$ for $s \in [0,t]$ by (2.17).
Since
$$\xi_y(t,y)= 1+t\bar u_x(y)+\int_0^t (t-s)\,\varphi_y(s,y)\,ds\,.
\eqno(2.23)$$
we find that in this case
$$\xi_y(t,y)=1+t\bar u_x(y)+{t^2 \over 4}\,\bar u_x^2(y)=
{1 \over 4}\,\Bigl(2+t\,\bar u_x(y)\Bigr)^2. \eqno(2.24)$$
In the remaining cases we have that $\bar u_x(y)=-\,{2 \over t_0}
\le -\,{2 \over t}$ for
some $t_0 \in (0,t]$. Therefore (2.17) yields $\varphi_y(s,y)
={1 \over 2}\,\bar u_x^2(y)$ for $s \in [0, t_0)$,
while $\varphi_y(s,y)=0$ for $s \in (t_0,t]$. From (2.23) we infer that
$$\xi_y(t,y)=1-{2t \over t_0}+{2t-t_0 \over t_0}=0.\eqno(2.25)$$
The relations (2.24)-(2.25) confirm the monotonicity of the map
$y \mapsto \xi(t,y)$. Since $\xi(0,x)=x$,
it remains to prove that $\lim_{y \to \pm\infty}\xi(t,y)=\pm\infty$
for any $t>0$. Fix $t>0$. Since
$\bar u_x \in L^2(\R)$, the Lebesgue measure $l(t)$ of the set
$\{y \in \R:\ \bar u_x(y) \le -\,{1 \over t}\}$
is finite. On the complement $C(t)$ of this set we obviously
have $\bar u_x(y) >-\,{1 \over t}$ and thus
$\xi_y(t,y) \ge {1 \over 4}$ by taking into account (2.24).
Therefore, given $x_2>x_1$, we infer that
$$\xi(t,x_2)-\xi(t,x_1)=\int_{x_1}^{x_2} \xi_y(t,y)\,dy \ge
\int_{[x_1,x_2]\,\cap\,C(t)} \xi_y(t,y)\,dy
\ge \int_{[x_1,x_2]\,\cap\,C(t)} {1 \over 4}\,dy\ge {x_2-x_1-l(t) \over 4}.$$
This proves the claim about the limiting behaviour of $\xi(t,\cdot)$
at $\pm\infty$. While for times $t$
up to the blow-up time $T_0$, given by (2.7), the map
$y \mapsto \xi(t,y)$ is an absolutely continuous
diffeomorphism of the real line, for $t \ge T_0$ this
map is nondecreasing and onto but is not necessarily
a bijection. Nevertheless, we would like to define the
solution $u$ by the formula (2.21) for all $t \ge 0$.

To show that $u$ is well-defined via (2.21), due to the
monotone and surjective character of the map
$y \mapsto \xi(t,y)$, it is sufficient to show that if
$\xi(t,y_1)=\xi(t,y_2)$ for some $y_2>y_1$, then
the values of $u$ given by (2.21) are also equal.
Indeed, we must have that $\xi(t,y)=\xi(t,y_1)$ for
all $y \in [y_1,y_2]$ and a glance at (2.24)-(2.25)
confirms that $\bar u_x(y) \le -\,{2 \over t}$ for
$y \in [y_1,y_2]$. This means that for every fixed
$y \in (y_1,y_2)$ we have $\bar u_x(y)=-\,{2 \over t_0(y)}$
for some $t_0(y) \in [0,t]$. Consequently $\varphi_y(s,y)
={1 \over 2}\,\bar u_x^2(y)\,\chi_{[0,t_0(y)]}(s)$ for
$s \in [0,t]$ and differentiation of the right-hand side of (2.21) yields
$$\partial_y\,\Bigl( \bar u(y)+ \int_0^t \varphi(s,y)\,ds \Bigr)
=\bar u_x(y) + \int_0^{t_0(y)} \,{1 \over 2}\, \bar u_x^2(y)\,ds
=-\,{2 \over t_0(y)}+{t_0(y) \over 2}\,{4 \over t_0^2(y)}=0$$
for $y \in (y_1,y_2)$. In particular, the values of the right-hand
side of (2.21) are equal when evaluated
at $(t,y_1)$ and at $(t,y_2)$. This proves that $u$ is well-defined.

The next step is to prove that for every $t \ge 0$, the map
$y \mapsto u(t,y)$ is continuous on $\R$ with distributional
derivative in $L^2(\R)$. Given $t \ge 0$ and
$y_0 \in \R$, let $I_0=\{ x \in \R:\ \xi(t,x)=y_0\}$.
The previously established
properties of the map $x \mapsto \xi(t,x)$ ensure that $I_0=[a,b]$
for some $a \le b$.
For any sequence $y_n \to y_0$, choose $x_n \in \R$ with $\xi(t,x_n)=y_n$. If
we show that $\min\,\{|x_n-a|,\,|x_n-b|\} \to 0$ as $n \to \infty$,
by the continuous
dependence on the $y$-variable of the right-hand side of (2.21), we infer
that
$$u(t,y_n)=u(t,\xi(t,n)) \to u(t,\xi(t,a))=u(t,\xi(t,b))=u(t,y_0)$$
since $\xi(t,x_n) \to \xi(t,a)=\xi(t,b)=y_0$. Thus $y \mapsto u(t,y)$
would be continuous
at $y_0$. If it would be possible that $\min\,\{|x_{n_k}-a|,\,|x_{n_k}-b|\}
\ge
\varepsilon >0$ for a sequence $n_k \to \infty$, then
$$|y_{n_k}-y_0|=|\xi(t,x_{n_k})-\xi(t,a)|=|\xi(t,x_{n_k})-\xi(t,b)|
\ge \min\,\{ y_0-\xi(t,a-\varepsilon),\,\xi(t,b+\varepsilon)\}>0$$
must hold by the definition of $[a,b]$ and the monotonicity property
of the function
$x \mapsto \xi(t,x)$. But this is a contradiction since $y_n \to y_0$
as $n \to \infty$.
We therefore proved the continuity of the map $y \mapsto u(t,y)$ for
every fixed
$y \in \R$. Actually, a glance at the previous considerations confirms
the continuity
of the map $u: \R_+ \times \R \to \R$ since $\xi: \R_+ \times \R \to \R$
is continuous
in view of (2.14). To show that for each $t \ge 0$ the distributional
derivative
$u_x(t,\cdot)$ belongs to $L^2(\R)$, due to the absolute continuity
of the nondecreasing surjective map $\xi(t,\cdot):\R \to \R$, we first
show that at every
point $y=\xi(t,x)$ where $\xi_x(t,x)>0$ exists, $u_x(t,y) \in \R$ exists.
Indeed,
at such a point $y$ the right-hand side of (2.21), formally equal to
$u_x(t,\xi(t,x))\,\cdot\,\xi_x(t,x)$, is differentiable with derivative
$$\bar u_x(x)+\int_0^t \varphi_y(t,x)\,ds=\bar u_x(x)
+{t \over 2}\,\bar  u_x^2(x),$$
since in view of (2.24)-(2.25) we must have $\bar u_x(x) > -\,{2 \over t}$ and
$\varphi_y(s,x)={1 \over 2}\,\bar u_x^2(x)$ for all $s \in [0,t]$. Since
$\xi_y(t,x)=1+t\,\bar u_x(x)+{t^2 \over 4}\,\bar u_x^2(x)$,
we infer that $u_x(t,y)$
exists, being given by the formula
$$u_x(t,y)={{\bar u_x(x)+{t \over 2}\,\bar u_x^2(x)}
\over {1+t\bar u_x(x)+{t^2 \over 4}\,\bar u_x^2(x)}}
={\bar u_x(x) \over {1+{t \over 2}\,\bar u_x(x)}},\eqno(2.26)$$
where $y=\xi(t,x)$. From (2.26) we deduce that for any interval $[x_1,x_2]$
where $\xi_x(t,x)>0$ a.e., we have
$$\eqalign{\int_{\xi(t,x_1)}^{\xi(t,x_2)} u_x^2(t,y)\,dy
&=\int_{\xi(t,x_1)}^{\xi(t,x_2)} u_x^2(t,\xi(t,x))\,\cdot\,\xi_x(t,x)\,dx \cr
&=\int_{y_1}^{y_2} {\bar u_x^2(x) \over {\Bigl(1+{t \over 2}
\,\bar u_x(x)\Bigr)^2}}\, \Bigl(1+t\,\bar u_x(x)+{t^2 \over 4}
\,\bar u_x^2(x)\Bigr)\,dx=\int_{y_1}^{y_2} \bar u_x^2(x)\,dx\cr}$$
if $y_1=\xi(t,x_1),\, y_2=\xi(t,x_2)$ and if we take into account (2.24).
Summing
up over such intervals, we obtain that
$$\int_{\R} u_x^2(t,x)\,dx = \int_{\{\bar u_x(x) >
 -\,{2 \over t}\}} \bar u_x^2(x)\,dx.\eqno(2.27)$$
In particular, the map $t \mapsto \| u_x(t,\cdot) \|_{L^2(\R)}$
is nondecreasing on $\R_+$. Moreover, we also have that
$$\int_{-\infty}^{\xi(t,y)} u_x^2(t,x)\,dx=\int_{\{ x \in (-\infty,y]:\
\bar u_x(x)>-\,{2 \over t}\}}\bar u_x^2(x)\,dx$$
and
$$\int_{\xi(t,y)}^\infty u_x^2(t,x)\,dx=\int_{\{ x \in [y,\infty):\
\bar u_x(x)>-\,{2 \over t}\}}\bar u_x^2(x)\,dx.$$
A comparison with (2.18) yields
$$\varphi(t,y)={1 \over 4} \int_{\R} \hbox{sign}\,\Bigl(\xi(t,y)-x\Bigr)\,
u_x^2(t,x)\,dx.\eqno(2.28)$$

Furthermore, if $\xi_x(t,x) >0$ exists, relation (2.20) ensures for
$y=\xi(t,x)$ the existence of $u_t(t,y)$, given by the formula
$$u_t(t,y)=-\,u_x(t,y)\,u(t,y)+\varphi(t,x)$$
obtained by differentiation and taking into acount (2.22).
In combination with (2.28),
this yields
$$\Bigl( u_t +uu_x\Bigr)(t,\xi(t,x))={1 \over 4}\,\int_{\R} \hbox{sign}
\Bigl(
\xi(t,x)-\zeta\Bigr)\,u_x^2(t,\zeta)\,d\zeta,$$
which is precisely (1.1) evaluated at $(t,\xi(t,x))$.
In view of the previously
established properties of the map $x \mapsto \xi(t,x)$
we deduce that the constructed function $u$ satisfies also
the condition (iii) of Definition 2.1.
Since the other properties required by Definition 2.1 were proved above, we
conclude that $u$ qualifies as a solution to (1.1) in the sense of
Definition 2.1. This
completes the proof of Theorem 1.$\hfill \diamondsuit$\bigskip

\vsk
\n{\medbf 3 - A distance functional}
\v
If $\bar u: \R \to \R$ is also bounded,
in addition to being continuous
and with distributional derivative $\bar u_x \in L^2(\R)$, then
the global solution $u(t,\cdot)$
constructed in Theorem 1 will be bounded at every fixed time $t \ge 0$.
More precisely, in view of (2.19) and (2.21) we have that
$$\sup_{t\geq 0,\,x \in \R}\,\big|u(t,x)\big| \leq \sup_{x \in \R}\,
\big|\bar u(x)
\big|
+{t \over 4}\,\int_{\R}
\bar u_x^2(x)\,dx\,.$$
Thus, if $\bar u : \R \to \R$ is a bounded continuous function with
distributional derivative $\bar u_x \in L^2(\R)$, then at each fixed time
$t \ge 0$,
the solution $u(t,\cdot)$ to (1.1), constructed in Theorem 1,
belongs to the
Banach space ${\cal X}$ of bounded continuous functions $f: \R \to \R$ with
distributional derivative $f_x \in L^2(\R)$, endowed with the norm
$$\| f \|_{\cal X}=\sup_{x \in \R}\,\{|f(x)|\}
+ \Bigl( \int_{\R} f_x^2(x)\,dx \Bigr)^{1 \over 2}.$$
The Banach space ${\cal X}$ seems suitable for (1.1) - see also [BZZ] where
a construction similar
to the one performed in Theorem 1 is presented. However, the map
$t \mapsto
u(t,\cdot)$ is generally not continuous from $\R_+$ to ${\cal X}$.
Indeed, if for some
$\tau >0$ we have that the set $\{x \in \R:\ \bar u_x(x)
=-\,{2 \over \tau}\}$ is of
positive Lebesgue measure, then a discontinuity occurs at time $t=\tau$ for
the map $t \mapsto u(t,\cdot) \in {\cal X}$ since from (2.27) we infer
that for $t < \tau$,
$$ \int_{\R} u_x^2(t,x)\,dx -\int_{\R} u_x^2(\tau,x)\,dx \ge
\int_{\{x \in x:\ \bar u_x(x)=-\,{2 \over \tau}\}} \bar u_x^2(x)\,dx>0.$$

Our aim will be to construct a distance functional in the space of
solutions to (1.1) with respect to which we will have both continuity with
respect to time as well as continuity with respect to the initial data for the
solutions to (1.1). More precisely, for non-smooth solutions the conservation law
(2.2) is replaced by
$$(v^2)_t+(uv^2)_x=-\mu\,,\eqno(3.1)$$
where $\mu$ is the positive measure on the $t$-$x$ plane defined as
$$\mu(\Omega)=\int_{\Big\{ \big(T(y)\,,~
\xi(T(y),y)\big)\,\in\,\Omega\Big\}}
\bar u^2_x(y)\,dy$$
for every open set $\Omega\subset \R_+\times\R$.
Here $T(y)$ is the blow-up time along the characteristic
curve starting at $y$, namely
$$T(y)\doteq \cases{ -2/\bar u_x(y)\qquad &if\quad $\bar
u_x(y)<0\,$,\cr
\infty\qquad &otherwise.\cr}
$$
For any $\bar u\in {\cal X}$, we can use the semigroup notation
$S_t\bar u\doteq u(t,\cdot)$ to denote the solution of (1.1)
constructed in Section 2.  Indeed
$$S_0\bar u=\bar u\,,\qquad S_{t+s}\bar u=S_t\big(S_s\bar u\big)\,.
\eqno(3.2)$$
To prove (3.2), we first show that
$$\xi_1(t+s,y)=\xi_2(t,\xi_1(s,y)),\qquad t,\,s \ge 0,\ y \in \R,\eqno(3.3)$$
where $\xi_2$ is the characteristic built upon the initial data $y \mapsto u(s,\xi_1(s,y))$. To check
(3.3), we view both expressions as functions of $t$. At $t=0$ they are both
equal to $\xi_1(s,y)$. For $t>0$, differentiation of (3.3) yields
$$\bar u(y)+\int_0^{t+s} \varphi_1(r,y)\,dr=u(s,\xi_1(s,y))+
\int_0^t \varphi_2(r,\xi_1(s,y))\,dr\eqno(3.4)$$
in view of (2.20). We use (2.21) to express the right-hand side of (3.4) as
$$\bar u(y)+ \int_0^s \varphi_1(r,y)\,dr+\int_0^t \varphi_2(r,\xi_1(s,y))\,dr.$$
Therefore, to get (3.4), which yields (3.3) by integration, it suffices to show that
$$\int_s^{t+s} \varphi_1(r,y)\,dr=\int_0^t \varphi_2(r,\xi_1(s,y))\,dr.\eqno(3.5)$$
To prove (3.5), we note that by (2.18),
$$\int_0^t \varphi_2(r,\xi_1(s,y))\,dr ={1 \over 4}\int_0^t \int_{\{x:\ u_x(s,x)>-{2 \over r}\}}\sgn\Big( \xi_1(s,y)-x\Big)\,u_x^2(s,x)\,dx\,dr$$
$$={1 \over 4}\int_0^t \int_{\{x:\ u_x(s,\xi_1(s,x))>-{2 \over r}\}}\sgn\Big( \xi_1(s,y)-\xi_1(s,x)\Big)\,u_x^2\Big(s,\xi_1(s,x)\Big)\,\partial_x\xi_1(s,x)\,dx\,dr$$
if we change variables $x \mapsto \xi_1(s,x)$. Taking now (2.9)-(2.10) into account, we
infer that
$$\eqalign{\int_0^t \varphi_2(r,\xi_1(s,y))\,dr &={1 \over 4}\int_0^t \int_{\{x:\ u_x(s,\xi_1(s,x))>-{2 \over r}\}}\sgn\Big( \xi_1(s,y)-\xi_1(s,x)\Big)\,\bar u_x^2(x)\,dx\,dr\cr
&={1 \over 4}\int_0^t \int_{\{x:\ u_x(s,\xi_1(s,x))>-{2 \over r}\}}\sgn(y-x)\,
\bar u_x^2(x)\,dx\,dr\cr}$$
since the function $x \mapsto \xi_1(s,x)$ is nondecreasing. But
$$u_x(s,\xi_1(s,x))={2\,\bar u_x(x) \over 2+s\,\bar u_x(x)}>-{2 \over r}\quad\hbox{if
and only if}\quad \bar u_x(x) >-\,{2 \over s+r}$$
since the function $y \mapsto \displaystyle{2y \over 2+ s\,y}$ is strictly increasing for
$y > -\,{2 \over s}$, so that in the end we get
$$\eqalign{\int_0^t \varphi_2(r,\xi_1(s,y))\,dr &={1 \over 4}\int_0^t \int_{\{x:\ \bar u_x(x)>-{2 \over r+s}\}}\sgn(y-x)\,\bar u_x^2(x)\,dx\,dr\cr
&={1 \over 4}\int_s^{t+s} \int_{\{x:\ \bar u_x(x)>-{2 \over \tau}\}}\sgn(y-x)\,
\bar u_x^2(x)\,dx\,d\tau\cr}\eqno(3.6)$$
where $\tau=r+s$. On the other hand, by (2.18),
$$\int_s^{t+s} \varphi_1(r,y)\,dr={1 \over 4}\int_s^{t+s} \int_{\{x:\ \bar u_x(x)>-{2 \over \tau}\}}\sgn(y-x)\,
\bar u_x^2(x)\,dx\,d\tau$$
so that (3.4) holds and (3.3) is proved. Knowing (3.3), to infer $S_{t+s}\bar u =S_t(S_s\bar u)$,
it suffices to show that
$$u(t+s,\xi_1(t+s,y))=u(t,\xi_2(t,\xi_1(s,y))).$$
But, by (2.21), the left-hand side is precisely
$$\bar u(y)+ \int_0^{t+s} \varphi(r,y)\,dr = \bar u(y)+ \int_0^s \varphi_1(r,y)\,dr
+ \int_s^{t+s}\varphi_1(r,y)\,dr
= u(s,\xi_1(s,y))+\int_s^{t+s}\varphi_1(r,y)\,dr,$$
which, taking into account (3.5), equals to
$$u(s,\xi_1(s,y))+\int_0^t \varphi_2(r,\xi_1(s,y))\,ds=u(t,\xi_2(t,\xi_1(s,y)))$$
in view of (2.21). This completes the proof of (3.2).

Notice that in general
the map $t\mapsto S_t\bar u$ is NOT continuous
from $[0,\infty[\,$ into ${\cal X}$.  It is thus interesting to identify some distance $J(u,v)$
which is well adapted to the evolution generated by
(1.1). More precisely, given an arbitrary constant $M$, in
this section we shall construct a functional $J(u,v)$
with the following property:
For any initial data $\bar u, \bar v\in {\cal X}$
with
$$\|\bar u_x\|_{\L^2}\leq M\,,\qquad\qquad
\|\bar v_x\|_{\L^2}\leq M\,,$$
the corresponding dissipative solutions $u,v$ constructed in Theorem 1
satisfy
$$J\big(u(t)\,,~v(t)\big)\leq e^{C_M\,t} J(\bar u,\,\bar v).$$

To begin the construction, consider
the metric space
$$X\doteq \Big(\R^2\times
\,]-\pi/2\,,~\pi/2]\Big)\cup\{\infty\}\eqno(3.7)$$
with distance
$$\eqalign{d\Big( (x,u,w), ~(\tilde x, \tilde u,\tilde w)\Big)
&\doteq \min\Big\{|x-\tilde x|+|u-\tilde u|+
\kappa_0\,|w-\tilde w|\,,\quad\kappa_0\,|\pi/2+w|+\kappa_0|\pi/2+\tilde
w|
\Big\}\,,\cr
d \Big( (x,u,w), ~\infty\Big)&=\kappa_0\,|\pi/2+w|
\,.\cr}\eqno(3.8)$$
Here $\kappa_0$ is a suitably large constant, whose precise value will
be specified later.
Notice that $X$ is obtained from the metric space $\R^2\times
[-\pi/2\,,~\pi/2]$ by identifying all points $(x,\,u,\, -\pi/2)$
into a single point, called ``$\infty$''.

Let $M(X)$ be the space of all bounded Radon measures on $X$.
To each function $u \in H^1_{loc}(\R)$ with $u_x \in L^2(\R)$
we now associate the measure
$\mu^u\in\M(X)$ defined as
$$\mu^u\big(\{\infty\}\big)=0\,,\qquad\qquad
\mu^u(A)=\int_{\big\{ x\in\R\,:~(x,\,u(x),\, \arctan \,u_x(x)\,)\in
A\big\}}
u_x^2(x)\,dx\eqno(3.9)$$
for every Borel set $A\subseteq\R^2\times \,]-\pi/2\,,~\pi/2]\,$.

As distance between two functions $u,v\in {\cal X}$ we now introduce a
kind of Kantorovich distance $J(u,v)$ related to an optimal transportation problem.
Call $\F$ the family of all triples $(\psi,\phi_1,\phi_2)$,
where $\phi_1,\phi_2:\R\mapsto [0,1]$ are simple Borel measurable maps
(that is, their range is a finite number of points and the preimage of each
such point is a Borel set) and $\psi:\R\mapsto \R$ is a nondecreasing absolute
continuous surjective map. Assuming that
$$\phi_1(x) \,u_x^2(x)=\psi'(x)\cdot\phi_2\big(\psi(x)\big)\,v_x^2
\big(\psi(x)\big)\qquad\qquad\hbox{for a.e.}~x\in\R\,,\eqno(3.10)$$
we define
$$\eqalign{J^{(\psi,\phi_1,\phi_2)}(u,v)\doteq &
\int d\Big( \big(x,\,u(x),\, \arctan u_x(x)\big)\,,~
\big(\psi(x),\,v(\psi(x)),\,\arctan v_x(\psi(x))\Big)
\cdot \phi_1(x)\,u_x^2(x)\,dx\cr
&+\int
d\Big( \big(x,\,u(x),\, \arctan u_x(x)\big)\,,~\infty\Big)\cdot
\big(1-\phi_1(x)\big)\, u_x^2(x)\,dx\cr
&
+\int  d\Big(
\big(\psi(x),\,v(\psi(x)\big),\,\arctan v_x(\psi(x)\big)\,,~\infty\Big)
\cdot \big(1-\phi_2(\psi(x))\big)\, v_x^2(\psi(x))\,\psi'(x)\,dx
\,.\cr}\eqno(3.11)$$
Observe that $(\psi,\phi_1,\phi_2)$
can be regarded as a {\bf transportation plan},
in order to transport the measure $\mu^u$
onto the measure $\mu^v$.
Since these two positive measures need not have the same total mass,
we allow some of the mass to be transferred to
the point $\infty$.  More precisely,
the mass transferred is $(1-\phi_1)\cdot\mu^u$
and $(1-\phi_2)\cdot\mu^v$.
The last two integrals in (3.11) account for the additional cost
of this transportation. Integrating (3.10) over the real line, one finds
$$\int_\R \phi_1(x) \,u_x^2(x)\,dx=
\int_\R \phi_2(y) \,v_x^2(y)\,dy\,.$$
We can thus transport the measure
$\phi_1 \mu^u$ onto $\phi_2\, \mu^v$ by a map $\Psi:\big(x,\, u(x)\,
\arctan u_x(x)\big)\mapsto
\big(y, \,v(y),\,\arctan v_x(y)\big)$,
with $y=\psi(x)$.
The associated cost is given by the first integral in (3.11). In this case the
measure $\phi_2\, \mu^v$ is obtained as the push-forward
of the measure $\phi_1 \mu^u$.
We recall that the {\bf push-forward} of a measure $\mu$ by a mapping
$\Psi$ is defined as $(\Psi\#\mu)(A)\doteq \mu(\Psi^{-1}(A))$
for every measurable set $A$.  Here $\Psi^{-1}(A)\doteq
\big\{ z\,:~\Psi(z)\in A\big\}$.

\v
We now define our distance functional
by optimizing over all transportation plans, namely
$$J(u,v)\doteq \inf_{(\psi,\phi_1,\phi_2)} \{J^{(\psi,\phi_1,\phi_2)}
(u,v)\}\,\eqno(3.12)$$
where the infimum is taken over all triples
$(\psi,\phi_1,\phi_2) \in\F$ such that (3.10) holds.
\v
To check that (3.12) actually defines a distance,
let $u,\,v,\,w\in {\cal X}$ be given
functions.
\v
\n{\bf 1.} Let us show that $J(u,v)=J(v,u)$. In order to do this, it is enough to prove that
for every triple $(\psi,\,\phi_1,\,\phi_2) \in \F$  satisfying (3.10) and every $\ve >0$, there is
a triple $(\eta,\varphi_1,\varphi_2) \in \F$ satisfying (3.10) such that $\eta: \R \to
\R$ is a strictly increasing absolutely continuous bijection and
$$\Big| J^{(\eta,\varphi_1,\varphi_2)}(u,v)-J^{(\psi,\phi_1,\phi_2)}(u,v)\Big|
 \le \ve.\eqno(3.13)$$
Indeed, given $(\psi,\phi_1,\phi_2) \in \F$
satisfying (3.10), define $\tilde{\psi}=\eta^{-1}$, $\tilde{\phi_1}=\varphi_2$, $\tilde{\phi_2}=
\varphi_1$. The properties of $\eta$ ensure the absolute continuity of $\tilde{\psi}$ (see [N])
so that we obtain $J^{(\tilde{\psi},\tilde{\phi_1},\tilde{\phi_2})}(v,u)=J^{(\eta,\varphi_1,\varphi_2)}(u,v)$
by performing the change of variables $x \mapsto \eta(x)$. Since $\ve>0$ was arbitrary,
we infer that $J(v,u) \le J(u,v)$. Interchanging the roles of $u$ and $v$ we get $J(u,v)=J(v,u)$.

To prove (3.13), it is convenient to view $\psi: \R \to \R$ as a maximal
monotone multifunction $\psi:\R\mapsto{\cal P}(\R)$ with domain and range $\R$.
Here ${\cal P}(\R)$ is the family of all subsets of $\R$. The conditions for a
multifunction $F:\R\mapsto{\cal P}(\R)$ to be maximal monotone
with domain and range $\R$ may be explicited as follows [Z]:

- for every $x \in \R$, the set $F(x) \subset \R$ is nonempty (i.e. the domain of $F$ is $\R$);

- for every $y \in \R$ there is at least some $x \in \R$ with $y \in F(x)$, expressing the fact that
the range of $F$ is $\R$;

- there are no couples $(x_1,y_1)$ and $(x_2,y_2)$ with $y_1 \in F(x_1)$ and
$y_2 \in F(x_2)$ such that $x_1<x_2$ and $y_2<y_1$, meaning that $F$ is monotone);

- if we associate to $F$ its graph $\{(x,y) \in \R^2:\ y \in F(x)\}$, then this graph
has no proper extension satisfying the first three properties (condition defining
the maximal monotonicity property).

We recall some important features presented by such maps [AA, Z]:

- the set $F(x)$ is an interval of the form $[a_x,b_x]$ with $a_x \le b_x$
for all $x \in \R$ and $a_x=b_x$ for all $x \in \R$, except perhaps an at most countable
set (so $F$ is singlevalued with the exception of at most countably many points);

- $F$ is a.e. differentiable, that is, for almost all $x_0 \in \R$ there exists $F'(x_0) \in \R$
such that
$$\lim_{x \to x_0,\, y \in F(x)} {y-F(x_0)-(x-x_0)\,F'(x_0) \over x-x_0}=0;$$

- we can define the inverse $F^{-1}: \R \to {\cal P}(\R)$ of $F$ by asking $y \in F^{-1}(x)$
if and only if $x \in F(x)$ and $F^{-1}$ is again a maximal monotone multifunction
with domain and range $\R$.

\noindent
Since the multifunction $\psi^{-1}$ is maximal montone, let $\{y_n\}$ be
the (at most countable) set of points where it is multivalued, that is, $\psi^{-1}(y)=[a_n,b_n]$
with $b_n>a_n$. Then $\psi(x)=y_n$ for $x \in [a_n,b_n]$ and, $\psi$ being absolutely
continuous, $\psi_x>0$ a.e. on $\R - \bigcup_{n}[a_n,b_n]$ since $\psi$ is strictly increasing
on this set. Given $\gamma>0$, the absolute continuity of $\psi: \R \to \R$ allows us to
choose some $\delta>0$ such that the total variation of $\psi$
over the union of disjoint closed intervals with the sum of their lengths less than $\delta$ is
less than $\gamma$ cf. [BGH]. On each interval $[a_n-{\delta \over 2^n},b_n+{\delta \over 2^n}]$
we replace $\psi$ with the linear function $\eta$ which takes the values
$\psi(a_n-{\delta \over 2^n})$, respectively $\psi(b_n+{\delta \over 2^n})$ at the endpoints.
By the way $a_n$ and $b_n$ were defined, we know that $\psi(b_n+{\delta \over 2^n})>
\psi(a_n-{\delta \over 2^n})$ so that $\eta'(x)$ is a positive constant on
$[a_n-{\delta \over 2^n},b_n+{\delta \over 2^n}]$ with
$$\sum_{n} \int_{a_n-{\delta \over 2^n}}^{b_n+{\delta \over 2^n}}\eta'(x)\,dx
\le \sum_{n} \Big(\psi(b_n+{\delta \over 2^n})-\psi(a_n-{\delta \over 2^n})\Big) \le \gamma.$$
Setting $\eta(x)=\psi(x)$ for $x \not \in
[a_n-{\delta \over 2^n},b_n+{\delta \over 2^n}]$, we obtain a strictly increasing
absolutely continuous bijection $\eta: \R \to \R$. Let us now show that the triple
$(\eta,\varphi_1,\varphi_2) \in \F$ satisfies both (3.10) and (3.13), where $\varphi_1,\,
\varphi_2$ are defined by setting $\varphi_1(x)=0$ for $x \in
[a_n-{\delta \over 2^n},b_n+{\delta \over 2^n}]$ and
$\varphi_1(x)=\phi_1(x)$  for $x \not \in
[a_n-{\delta \over 2^n},b_n+{\delta \over 2^n}]$, while
$\varphi_2(\eta(x))=\phi_2(\psi(x)))$ for $x \not \in
[a_n-{\delta \over 2^n},b_n+{\delta \over 2^n}]$ and
$\varphi_2(\eta(x))=0$  for $x \in
[a_n-{\delta \over 2^n},b_n+{\delta \over 2^n}]$. On the complement
of the set $\bigcup_n [a_n-{\delta \over 2^n},b_n+{\delta \over 2^n}]$
relation $(3.10)$ clearly holds a.e. for $(\eta,\varphi_1,\varphi_2)$, being unmodified
from (3.10) for $(\psi,\phi_1,\phi_2)$. If
$x \in [a_n-{\delta \over 2^n},b_n+{\delta \over 2^n}]$, then (3.10) for
$(\eta,\varphi_1,\varphi_2)$ holds again since both
sides are zero as $\varphi_1(x)=\varphi_2(\eta(x))=0$ in this case. Finally, to
check (3.13), notice that if we denote
$$E_\delta=\bigcup_{n} \Big\{[a_n-{\delta \over 2^n},a_n]\,\cup\,[b_n,b_n+{\delta \over 2^n}]\Big\},\qquad A=\bigcup_{n}\,[a_n,b_n],$$
then
$$\Big| J^{(\eta,\varphi_1,\varphi_2)}(u,v)-J^{(\psi,\phi_1,\phi_2)}(u,v)\Big|
 \le 2\kappa_0\pi  \int_{E_\delta} u_x^2\,dx +2\kappa_0\pi  \int_{E_\delta\,\cup\,A} v_x^2(\eta(x))\,\eta'(x)\,dx.\eqno(3.14)$$
Indeed, the distance $d$ is less that $2\kappa_0\pi$ and the integrands in $J^{(\eta,\varphi_1,\varphi_2)}(u,v)$ and
$J^{(\psi,\phi_1,\phi_2)}(u,v)$ agree on the complement
of the set $\bigcup_n [a_n-{\delta \over 2^n},b_n+{\delta \over 2^n}]$ by definition.
Also, for a.e. $x \in [a_n,b_n]$ we have $\phi_1(x)\,u_x^2=0$ by
(3.10) as $\psi'(x)=0$, and $\varphi_1(x)=0$ by its definition. We obtain (3.14).
Since the absolutely continuous map $\eta$ maps $E_\delta \,\cup\,A$ into a set
of Lebesgue measure less than $\gamma$, and $u_x^2,\,v_x^2 \in L^1(\R)$, from
(3.14) we infer (3.13) by choosing $\delta >0$ and $\gamma>0$ small enough. This
completes the argumentation needed to show that $J(u,v)=J(v,u)$.
\v
\n{\bf 2.} Choosing $\psi(x)=x$, $\phi_1(x)=\phi_2(x)=1$, we immediately see that
$J(u,u)=0$. Moreover, we have $J(u,v)>0$ if $u \not = v$. To check this,
note that $J(u,v)=0$ implies that there is a sequence $(\psi^n,\phi_1^n,\phi_2^n)$
along which $J^{(\psi^n,\phi_1^n,\phi_2^n)}(u,v) \to 0$. The second term in (3.11)
yields
$$\Big( {\pi \over 2}+\arctan u_x(x) \Big)\,\Big( 1 - \phi_1^n(x)\Big)\,u_x^2(x) \to 0
\quad \hbox{in}\quad L^1(\R),$$
so that along a subsequence $(1-\phi_1^{n_k})\,u_x^2 \to 0$ a.e. on $\R$ since
$u_x>-\infty$ a.e. On the set $S=\{x \in \R:\ u_x(x) \not =0\}$ we therefore have $\phi_1^{n_k} \to 1$ a.e. Moreover, the first term in (3.11) forces
$$\phi_1^{n_k}(x)u_x^2(x) \cdot\, \min\,\Big\{ |x-\psi^{n_k}(x)| +
|u(x)-v(\psi^{n_k}(x))|+\kappa_0 |\arctan u_x(x) - \arctan v_x(\psi^{n_k}(x))|,\ $$
$$\kappa_0  \,\Big[{\pi \over 2} + \arctan u_x(x)+
{\pi \over 2} + \arctan v_x(\psi^{n_k}(x)) \Big]\Big\} \to 0 \quad \hbox{in} \quad
L^1(\R). \eqno(3.15)$$
Since $u_x>-\infty$ a.e. ensures
$${\pi \over 2} + \arctan u_x(x)+
{\pi \over 2} + \arctan v_x(\psi^{n_k}(x)) \ge {\pi \over 2} + \arctan u_x(x)>0\quad
\hbox{a.e. on}\quad \R,$$
we infer from (3.15), by passing to another subsequence, that
$$|x-\psi^{n_k}(x)| +
|u(x)-v(\psi^{n_k}(x))| \to 0 \quad\hbox{a.e. on}\quad S.$$
In view of the continuity of $v$, $\psi^{n_k}(x) \to x$ a.e. on $S$ guarantees $v(\psi^{n_k}(x)) \to v(x)$ a.e. on $S$ so that $u=v$ a.e. on $S$ since also
$v(\psi^{n_k}(x)) \to u(x)$ a.e. on $S$. Repeating this argument with the roles
of $u$ and $v$ reversed, we find that $u=v$ a.e. on the set $\{x \in \R:\ v_x \not =0\}$.
Combining this with the previous conclusion, we have $u=v$ a.e. on the complement
of the set $\{x \in \R: u_x=v_x=0\}$. Since $u_x,\,v_x \in L^2(\R)$, this is possible
only if $u=v$ on $\R$. Thus $J(u,v)=0$ if and
only if $u=v$.
\v
\n{\bf 3.} Finally, to prove the triangle inequality, it suffices to show that for every
choice of $(\psi^\flat,\phi_1^\flat,\phi_2^\flat)$ satisfying (3.10), and of $(\psi^\sharp,\phi_1^\sharp,\phi_2^\sharp)$ satisfying (3.10) for $(v,w)$, the
triplet $(\psi,\phi_1,\phi_2)$ defined by
$$\psi(x)=\psi^\sharp( \psi^\flat(x)),\qquad \phi_1(x)=\phi_1^\flat(x)\cdot \phi_1^\sharp( \psi^\flat(x)),\quad \phi_2(y)=\phi_2^\sharp(y) \cdot \phi_2^\flat( \psi^\flat(x)),$$
satisfies (3.10) for $(u,w)$ and
$$J^{(\psi,\phi_1,\phi_2)}(u,w) \le J^{(\psi^\flat,\phi_1^\flat,\phi_2^\flat)}(u,v)
+ J^{(\psi^\sharp,\phi_1^\sharp,\phi_2^\sharp)}(v,w).\eqno(3.16)$$

Notice that composing the relation (3.10) for $(v,w)$ a.e. to the right with $\psi^\flat$,
and multiplying the outcome by $\phi_2^\flat \circ \psi^\flat \cdot (\psi^\flat)'$,
we infer that (3.10) holds a.e. on $\R$ for $(u,w)$ with our choice of $(\psi,\phi_1,\phi_2)$
and we can now concentrate on proving (3.16).

To simplify matters, we introduce the following notation
$$P_1=(x,\,u,\,\arctan u_x),\quad  P_2 =(\psi^\flat,\,v \circ \psi^\flat,\,
\arctan v_x \circ \psi^\flat),\quad  P_3=(\psi,\,w \circ \psi,\, \arctan w_x \circ \psi),$$
$$m_1=u_x^2,\qquad  m_2=v_x^2 \circ \psi^\flat \cdot (\psi^\flat)',
\qquad  m_3=w_x^2 \circ \psi \cdot \psi'.$$
The relations of type (3.10) yield then that a.e. on $\R$,
$$\phi_1^\flat \cdot m_1=\phi_2^\flat \circ \psi^\flat \cdot m_2,\qquad
\phi_1^\sharp \circ \psi^\flat \cdot m_2=\phi_2^\sharp \circ \psi \cdot m_3,
\qquad \phi_1\cdot m_1=\phi_2\circ \psi \cdot m_3.\eqno(3.17)$$
Also,
$$J^{(\psi,\phi_1,\phi_2)}(u,w) =\int_{\R} \Big\{ d(P_1,P_3)\cdot \phi_1 m_1 + d(P_1,\infty)\cdot (1- \phi_1) m_1+d(P_3,\infty)\cdot (1-\phi_2 \circ \psi)m_3\Big\}\,dx,$$
$$J^{(\psi^\flat,\phi_1^\flat,\phi_2^\flat)}(u,v) =\int_{\R} \Big\{ d(P_1,P_2)\cdot \phi_1^\flat m_1 + d(P_1,\infty)\cdot (1- \phi_1^\flat) m_1+d(P_2,\infty)\cdot (1-\phi_2^\flat \circ \psi^\flat)m_2\Big\}\,dx,$$
$$J^{(\psi^\sharp,\phi_1^\sharp,\phi_2^\sharp)}(v,w) = \int_{\R} \Big\{ d(P_2,P_3)\cdot
\phi_1^\sharp \circ \psi^b m_1 + d(P_2,\infty)\cdot (1- \phi_1^\sharp \circ \psi^\flat) m_2+d(P_3,\infty)\cdot (1-\phi_2^\sharp \circ \psi)m_3\Big\}\,dx,$$
the last relation being obtained after the change of variables
$x \mapsto \psi^\flat(x)$ in the integral.
We will prove (3.16) by deriving an appropriate inequality valid a.e. pointwise between
the integrands in the previous expressions. Since
$$(1- \phi_2^\flat \circ \psi^\flat)(1-\phi_1^\sharp \circ \psi^\flat) \ge 0,$$
we have
$$1- \phi_2^\flat \circ \psi^\flat+1-\phi_1^\sharp \circ \psi^\flat \ge
\phi_2^\flat \circ \psi^\flat(1-\phi_1^\sharp \circ \psi^\flat) +\phi_1^\sharp \circ \psi^\flat (1-\phi_2^\flat \circ \psi^\flat).$$
Multiplication of both sides by $d(P_2,\infty)\cdot m_2$ leads to
$$d(P_2,\infty)\cdot (1- \phi_2^\flat \circ \psi^\flat)m_2+d(P_2,\infty)\cdot (1-\phi_1^\sharp \circ \psi^\flat)m_2 \ge d(P_2,\infty)\cdot
\phi_1^\flat (1-\phi_1^\sharp \circ \psi^\flat)m_1 $$
$$\qquad +\ d(P_2,\infty)\cdot
\phi_1^\sharp \circ \psi^\flat (1-\phi_2^\flat \circ \psi^\flat)m_2\eqno(3.18)$$
in view of (3.17). Multiply now the inequalities
$$d(P_1,P_2)-d(P_1,\infty)+d(P_2,\infty) \ge 0,\qquad d(P_2,P_3) - d(P_3,\infty)+d(P_2,\infty) \ge 0,$$
by $\phi_1^\flat(1-\phi_1^\sharp \circ \psi^\flat)m_1$, respectively
$\phi_1^\sharp \circ \psi^\flat(1-\phi_2^\flat \circ \psi^\flat)m_2$, and add them up.
The outcome yields in combination with (3.18) that
$$d(P_1,P_2)\cdot \phi_1^\flat(1-\phi_1^\sharp \circ \psi^\flat)m_1-
d(P_1,\infty)\cdot \phi_1^\flat(1-\phi_1^\sharp \circ \psi^\flat)m_1
+d(P_2,P_3)\cdot \phi_1^\sharp \circ \psi^\flat(1-\phi_2^\flat \circ \psi^\flat)m_2 $$
$$\quad + d(P_2,\infty)\cdot (1- \phi_2^\flat \circ \psi^\flat)m_2+
d(P_2,\infty)\cdot (1-\phi_1^\sharp \circ \psi^\flat)m_2 \ge
d(P_3,\infty)\cdot \phi_1^\sharp \circ \psi^\flat(1-\phi_2^\flat \circ \psi^\flat)m_2 .$$
Adding to both sides the quantity
$$d(P_1,\infty)\cdot m_1+d(P_3,\infty)\cdot m_3 + d(P_1,P_2)\cdot
\phi_1^\flat \cdot\phi_1^\sharp \circ \psi^\flat \cdot m_1 -
d(P_1,\infty)\cdot \phi_1^\flat\cdot \phi_1^\sharp \circ \psi^\flat\cdot m_1$$
$$\qquad +d(P_2,P_3)\cdot \phi_1^\sharp \circ \psi^\flat \cdot
\phi_2^\flat \circ \psi^\flat \cdot m_2 - d(P_3,\infty)\cdot
\phi_2^\sharp \circ \psi \cdot m_3$$
we deduce by (3.17) that the integrand of $J^{(\psi^\flat,\phi_1^\flat,\phi_2^\flat)}(u,v)
+ J^{(\psi^\sharp,\phi_1^\sharp,\phi_2^\sharp)}(v,w)$, equal a.e. precisely to the
left-hand side of the new inequality, is a.e. pointwise larger than
$$d(P_3,\infty)\cdot \phi_1^\sharp \circ \psi^\flat(1-\phi_2^\flat \circ \psi^\flat)m_2+
d(P_1,\infty)\cdot m_1+d(P_3,\infty)\cdot m_3 + d(P_1,P_2)\cdot
\phi_1^\flat \cdot\phi_1^\sharp \circ \psi^\flat \cdot m_1$$
$$ \qquad -d(P_1,\infty)\cdot \phi_1^\flat\cdot \phi_1^\sharp \circ \psi^\flat\cdot m_1 +d(P_2,P_3)\cdot \phi_1^\sharp \circ \psi^\flat \cdot
\phi_2^\flat \circ \psi^\flat \cdot m_2 - d(P_3,\infty)\cdot
\phi_2^\sharp \circ \psi \cdot m_3.$$
Taking into account (3.17) and the definition $\phi_1=\phi_1^\flat
\cdot \phi_1^\sharp \circ \psi^\flat$, we see that the above expression equals
$$d(P_1,\infty)\cdot (1-\phi_1)m_1+d(P_3,\infty)\cdot (1-\phi_2 \circ \psi)m_3
+ \Big( d(P_1,P_2)+d(P_2,P_3)\cdot \phi_1 m_1\Bigr)$$
$$\qquad \ge d(P_1,\infty)\cdot (1-\phi_1)m_1+d(P_3,\infty)\cdot (1-\phi_2 \circ \psi)m_3 +d(P_1,P_3)\cdot\phi_1m_1.$$
The lower estimate is a.e. precisely the integrand in $J^{(\psi,\phi_1,\phi_2)}(u,w)$ and (3.16) holds. The proof that $J$ satisfies the triangle inequality is therefore completed.

\vs
In the remainder of this section we examine
how the distance $J(\cdot,\cdot)$
behaves in connection with solutions
of the equation (1.1).

\vs
\n{\bf Continuity w.r.t.~time.}
Let $t\mapsto u(t)$ be the solution of (1.1) constructed in Section 2.
For any fixed $t>0$, we define a
transportation plan of $\mu^{\bar u}$ to $\mu^{u(t)}$ by setting
$$\psi(x)\doteq \xi(t,x)\,,\qquad
\phi_1(x)\doteq\cases{ 1\quad &if\quad $T(x)>t\,$,\cr
 0\quad &if\quad $T(x)\leq t\,$,\cr}\qquad
\phi_2(x)\equiv 1\,.\eqno(3.19)$$
Relation (3.6) follows from (2.9)-(2.10) on $\{T(x)>t\}$ and from (2.25) on $\{T(x) \le t\}$.
The cost of this plan is estimated by
$$\eqalign{J^{(\psi,\phi_1,\phi_2)}\big(\bar u\,,~u(t)\big)&\leq
\int_{\big\{ T(x)>t\big\}}
\Big\{ \big|x-\xi(t,x)\big|+ \big|\bar u(x)-u(t,\xi(t,x))\big|\cr
&\qquad\qquad\qquad
+\kappa_0 \big|\arctan \bar u_x(x)-\arctan  u_x(t,\xi(t,x))
\big|\Big\}\, \bar u_x^2(x)\,dx\cr
&\qquad +\int_{\big\{ T(x)\leq t\big\}}
\big| \pi/2 +\arctan \bar u_x(x)\big|\,\bar u_x^2(x)\,dx\,.\cr} \eqno(3.20)$$
By (2.4) we have that a.e.
$$\left|{d\over dt} \arctan u_x\big(t,\xi(t,x)\big)\right|
=\left| { {d \over dt} u_x\big(t,\xi(t,x)\big) \over 1+u_x^2\big(t,\xi(t,x)\big) }\right|\le {1 \over 2}.
\eqno(3.21)$$
An integration on $[0,t]$ yields
$$\left|\,\arctan \bar u_x(x)-\arctan u_x\big(t,\xi(t,x)\big)\right| \le {t \over 2},\qquad t \ge 0.\eqno(3.22)$$
On the other hand, using (2.20), we get
$$ |x - \xi(t,x)|  \le t\,|\bar u(x)|+ \int_0^t (t-s)\,|\varphi(s,x)|\,ds \le t\,|\bar u(x)|+{t^2 \over 8} \int_{\R} \bar u_x^2(x)\,dx,\qquad t \ge 0,\ x \in \R, \eqno(3.23)$$
if we take into account (2.19). From (2.21) and (2.19), we also infer
$$\left| \bar u(x) - u\big(t,\xi(t,x)\big) \right| \le \int_0^t |\varphi(s,y)|\,dy \le
{t \over 4}\, \int_{\R} \bar u_x^2(x)\,dx,\qquad t \ge 0,\ x \in \R.\eqno(3.24)$$
To estimate the last term in (3.20), notice that
$$\{x \in \R:\  T(x) \le t\}=\{ x\in \R:\ -\,\displaystyle{2 \over \bar u_x(x)} \le t\}
=\{ x \in \R:\ \bar u_x(x) \le -\,\displaystyle{2 \over t}\},\qquad t >0.\eqno(3.25)$$
Furthermore, since $\displaystyle\lim_{x \to - \infty}x({\pi \over 2}+\arctan x)=-1$, there is
a constant $c>0$ such that
$$0 \le {\pi \over 2} + \arctan y \le {c \over |y|},\qquad y \le -1,$$
whereas
$$\left| {\pi \over 2} + \arctan y \right| \,y^2 \le \pi\qquad\hbox{if}\quad -1 \le y \le 0,$$
so that
$$\left| {\pi \over 2}+\arctan \bar u_x(x) \right|\, \bar u_x^2(x) \le \pi + c \,|\bar u_x(x)|\qquad
\hbox{if}\quad \bar u_x(x) \le -\,{2 \over t}.\eqno(3.26)$$
On the other hand, if $\bar u_x(x) \le -\,\displaystyle{2 \over t}$, then $t^2 \bar u_x^2(x) \ge 4$
so that
$$\int_{\{T(x) \le t\}} 1\,dx \le {t^2 \over 4}\,\int_{\{T(x) \le t\}} \bar u_x^2(x)\,dx.\eqno(3.27)$$
From (3.25)-(3.27) we infer that
$$\eqalign{&\int_{\{T(x) \le t\}} \left| {\pi \over 2}+\arctan \bar u_x(x) \right|\, \bar u_x^2(x)\,dx
\le {\pi t^2 \over 4} \, \Vert \bar u_x \Vert_{L^2}^2 +c \int_{\{T(x) \le t\}} |\bar u_x(x)|\,dx \cr
&\le  {\pi t^2 \over 4} \, \Vert \bar u_x \Vert_{L^2}^2 +c \Big(\int_{\{T(x) \le t\}} 1\,dx\Big)^{1 \over 2} \Big(\int_{\{T(x) \le t\}} \bar u_x^2(x)\,dx\Big)^{1 \over 2}
\le {\pi t^2 \over 4} \, \Vert \bar u_x \Vert_{L^2}^2 +
{ct \over 2}\, \Vert \bar u_x \Vert_{L^2}^2.}$$
By (3.20), (3.22)-(3.25) and the previous inequality we conclude
$$J^{(\psi,\phi_1,\phi_2)}\big(\bar u\,,~u(t)\big)\le \Big({\pi t \over 4} +
 {c +\kappa_0 \over 2}+ \Vert \bar u \Vert_{L^\infty} +{t +2 \over 8}\,\Vert \bar u_x \Vert_{L^2}^2\Big)\,t\,\Vert \bar u_x \Vert_{L^2}^2,\qquad t \ge 0.\eqno(3.28)$$
It is now clear that each semigroup trajectory
$t\mapsto S_t\bar u$ is Lipschitz continuous as a map from
$[0,\infty[\,$ into the metric space $X$ equipped with
our distance functional $J$.
The Lipschitz constant remains uniformly
bounded as $\bar u$ ranges over bounded subsets of $X$.

\vs
\n{\bf Continuity w.r.t.~the initial data.}
We now consider two distinct solutions
and study how the distance $J\big(u(t)\,,~\tilde u(t)\big)$
varies in time.
Recall that the solution $u=u(t,x)$
is computed by (2.20)--(2.22), also in the case where
the gradient blows up.  The same formula of course holds
for $\tilde u$.
Let $(\psi_0,\phi_{1,0},\phi_{2,0})$ be an optimal transportation plan
of the measure $\mu^{u(0)}$ to the measure $\mu^{\tilde u(0)}$. In
view of the approximation property established in (3.13), we can restrict
our attention to the case when $\psi_0$ is strictly increasing on $\R$.
For any $t>0$, we define a
transportation plan $(\psi^t,\phi_1^t,\phi_2^t)$
of the measure $\mu^{u(t)}$ to  $\mu^{\tilde u(t)}$ as follows:
$$\psi^t\big(\xi(t,y)\big)\doteq \tilde\xi(t,~\tilde y)\quad\hbox{for}\quad \tilde y =\psi_0(y),$$
$$\phi_1^t\big(\xi(t,y)\big)\doteq \cases{ \phi_{1,0}(y)\qquad
&if\quad $T(y)>t\quad\hbox{and}\quad~\widetilde T(\tilde y)>t\quad\hbox{for}\ \tilde y = \psi_0(y)\,$,\cr
0\qquad &otherwise,\cr}$$
$$\phi_2^t\big(\tilde\xi\big(t,~\tilde y)\big)
\doteq \cases{ \phi_{2,0}(\tilde y)\qquad
&if\quad $T(y)>t\quad\hbox{and}\quad~\widetilde T(\tilde y)>t\quad\hbox{for}\ y = \psi_0^{-1}(\tilde y)\,$,\cr
0\qquad &otherwise.\cr}$$
If initially the point $y$ is mapped to $\tilde y= \psi_0(y)$,
then at any later time $t>0$ the point $\xi(t,y)$ along the $u$-characteristic
starting from $y$ is sent to the point $\tilde\xi(t,\,\tilde y)$
along the $\tilde u$-characteristic starting from $\tilde y=\psi_0(y)$.
We thus transport the mass from the point
$\Big(\xi(t,y)\,,~u\big(t,\xi(t,y)\big)\,,~\arctan u_x\big(t,\xi(t,y)\big)\Big)$
to the corresponding point $\Big(\tilde \xi(t,\tilde y)\,, ~\tilde u\big(t,\tilde \xi(t,\tilde
y)\big)\,,
~\arctan \tilde u_x\big(t,\tilde \xi(t,\tilde y)\big)\Big)\quad\hbox{with}\quad \tilde y= \psi_0(y),$ except in the case where blow up has occurred within time $t$
along one (or both) characteristics
$\xi(\cdot,y),\,\tilde\xi(\cdot,\tilde y)$.
In this later case, the mass is transported to the point $\infty$.

To check (3.10), it suffices to show that a.e.
$$\phi_1^t(\xi(t,y))\cdot u_x^2(t,\xi(t,y))\cdot \xi_x(t,y)=\phi_2^t\Big(\psi^t(\xi(t,y))\Big)\cdot
(\psi^t)'(\xi(t,y)) \cdot \xi_x(t,y) \cdot u_x^2\Big(t,\psi^t(\xi(t,y))\Big).$$
Since the relations $\tilde y =\psi_0(y)$, $\psi^t(\xi(t,y))=\tilde \xi(t,\,\tilde y)$, and
$(\psi^t)'(\xi(t,y)) \cdot \xi_x(t,y)=\tilde \xi_x(t,\,\psi_0(y))\cdot \psi_0'(y)$ all hold a.e., the
desired identity holds a.e. on the complement of the set $\{y:\ \tilde y= \psi_0(y),\ T(y)>t,\ \tilde
T(\tilde y)>t\}$ where both sides equal zero since $\phi_1^t(\xi(t,y))=\phi_2^t(\tilde \xi(t,\tilde y))=0$. The identity holds also a.e. on the set $\{y:\ \tilde y= \psi_0(y),\ T(y)>t,\ \tilde
T(\tilde y)>t\}$ since there, in view of (2.9)-(2.10), it practically amounts to relation (3.10)
for $(\phi_{1,0},\,\phi_{2,0},\,\psi_0)$.

In the following, our main goal is to provide an estimate on the
time derivative of the function
$$J(t)=J^{(\psi^t,\phi_1^t,\phi_2^t)} \big(u(t),\tilde u(t)\big).$$
Throughout the remainder of this section, by $\{\tilde T(\tilde y)^>_\le t\}$ we understand the set
of all $y \in \R$ such that $\psi_0(y)=\{\tilde y\}$ and $\tilde T(\tilde y)^>_\le t$. Since $u_x^2(t,\xi(t,y))\cdot \xi_x(t,y)=\bar u_x^2(y)$ on $\{T(y)>t\}$ by (2.9)-(2.10)
and $\tilde u_x^2(t,\,\tilde \xi(t, \tilde y))\cdot \tilde \xi_x(t,\tilde y)=\tilde u_x^2(0,\tilde y)$
on $\{\tilde T(\tilde y)>t\}$, while
$\psi^t(\xi(t,y))=\tilde\xi(t,\tilde y)$ for $\psi_0(y)=\{\tilde y\}$,
performing the change of variables $y \mapsto \xi(t,y)$, we see that
$$\eqalign{J(t) &= \int_{\{T(y)>t,\,\tilde T(\tilde y)>t\}}  \min\,\Big\{ |\xi(t,y)-\tilde \xi(t,\tilde y)|
+ |u(t,\xi(t,y))-\tilde u(t,\tilde \xi(t,\tilde y))|\cr
&\qquad\qquad + \kappa_0 |\arctan u_x(t,\xi(t,y)) - \arctan \tilde u_x(t,\tilde \xi(t,\tilde y))|,\cr
&\qquad\qquad \kappa_0 \Big( \pi + \arctan u_x(t,\xi(t,y)) + \arctan \tilde u_x(t,\tilde \xi(t,\tilde y))\Big)\Big\} \,\phi_{1,0}(y)\,\bar u_x^2(y)\,dy\cr
&+\kappa_0 \int_{\{ T(y) >t,\, \tilde T(\tilde y)>t\}} \Big( {\pi \over 2} + \arctan u_x(t,\xi(t,y))\Big)\,
\Big(1-\phi_{1,0}(y)\Big)\,\bar u_x^2(y)\,dy\cr
&+\kappa_0 \int_{\{ T(y) \le t \ \hbox{or}\  \tilde T(\tilde y) \le t\}}
\Big( {\pi \over 2} + \arctan u_x(t,\xi(t,y))\Big)\,\bar u_x^2(y)\,dy\cr
&+\kappa_0 \int_{\{ T(y) >t,\, \tilde T(\tilde y)>t\}} \Big( {\pi \over 2} + \arctan \tilde u_x(t,\tilde \xi(t,\,\tilde y))\Big)\,
\Big(1-\phi_{2,0}(\tilde y)\Big)\,\tilde u_x^2(0,\tilde y)\,\psi_0'(y)\,dy\cr
&+\kappa_0 \int_{\{ T(y) \le t \ \hbox{or}\  \tilde T(\tilde y) \le t\}}
\Big( {\pi \over 2} + \arctan \tilde u_x(t,\,\tilde \xi(t,\,\tilde y))\Big)\,
\tilde u_x^2(0,\tilde y)\,\psi_0'(y)\,dy.\cr}$$
To simplify notation, let
$$S(t)=\{ T(y)>t,\, \tilde T(\tilde y)>t\},\qquad S^c(t)=\R -S(t),\eqno(3.29)$$
$$E(t,y)=\min\,\Big\{ |u(t,\xi(t,y))-\tilde u(t,\tilde \xi(t,\tilde y))|+ \kappa_0 |\arctan u_x(t,\xi(t,y))  - \arctan \tilde u_x(t,\tilde \xi(t,\tilde y))|$$
$$+  |\xi(t,y)-\tilde \xi(t,\tilde y)|
,\quad \kappa_0 \Big( \pi +
\arctan u_x(t,\xi(t,y)) + \arctan \tilde u_x(t,\tilde \xi(t,\tilde y))\Big)\Big\}. \eqno(3.30)$$
Since for $h \ge 0$ we have
$$S(t+h) \subset S(t),\qquad S^c(t) \subset S^c(t+h),\eqno(3.31)$$
we deduce that
$$J(t+h)-J(t) = \int_{S(t)} \Big( E(t+h,y)-E(t,y)\Big) \,\phi_{1,0}(y)\,\bar u_x^2(y)\,dy$$
$$- \int_{S(t) \setminus S(t+h)} E(t+h,y)\,\phi_{1,0}(y)\,\bar u_x^2(y)\,dy$$
$$ + \kappa_0 \int_{S(t)} \Big( \arctan u_x(t+h,\xi(t+h,y))-\arctan u_x(t,\xi(t,y))\Big)
\,\Big(1-\phi_{1,0}(y)\Big)\,\bar u_x^2(y)\,dy$$
$$ +  \kappa_0 \int_{S(t)} \Big( \arctan \tilde u_x (t+h,\tilde \xi(t+h,\tilde y))-
\arctan \tilde u_x (t,\tilde \xi(t,\tilde y))\Big)
\,\Big(1-\phi_{2,0}(\tilde y)\Big)\, \tilde u_x^2(0,\tilde y)\,\psi_0'(y)\,dy$$
$$ - \kappa_0 \int_{S(t)\setminus S(t+h)} \Big( {\pi \over 2} + \arctan u_x(t+h,\xi(t+h,y))\Big)
\,\Big(1-\phi_{1,0}(y)\Big)\,\bar u_x^2(y)\,dy$$
$$-\kappa_0 \int_{S(t) \setminus S(t+h)}  \Big( {\pi \over 2} +
\arctan \tilde u_x (t+h,\tilde \xi(t+h,\tilde y))\Big) \,\,\Big(1-\phi_{2,0}(\tilde y)\Big)\, \tilde u_x^2(0,\tilde y)\,\psi_0'(y)\,dy$$
$$ + \kappa_0 \int_{S^c(t)} \Big( \arctan u_x(t+h,\xi(t+h,y))-\arctan u_x(t,\xi(t,y))\Big)
\,\bar u_x^2(y)\,dy$$
$$ +  \kappa_0 \int_{S^c(t)} \Big( \arctan \tilde u_x (t+h,\tilde \xi(t+h,\tilde y))-
\arctan \tilde u_x (t,\tilde \xi(t,\tilde y))\Big)
\, \tilde u_x^2(0,\tilde y)\,\psi_0'(y)\,dy$$
$$ + \kappa_0 \int_{S^c(t+h)\setminus S^c(t)} \Big( {\pi \over 2} + \arctan u_x(t+h,\xi(t+h,y))\Big)
\,\bar u_x^2(y)\,dy$$
$$+\kappa_0 \int_{S^c(t+h) \setminus S^c(t)}  \Big( {\pi \over 2} +
\arctan \tilde u_x (t+h,\tilde \xi(t+h,\tilde y))\Big) \,\tilde u_x^2(0,\tilde y)\,\psi_0'(y)\,dy.
\eqno(3.32)$$
Noticing that $S(t) \setminus S(t+h) =S(t)\, \cap \,S^c(t+h) =S^c(t+h) \setminus S^c(t)$,
we see that the combination of the fifth and ninth terms above, with that of the sixth and tenth,
added to the second term, amount to
$$\kappa_0 \int_{S(t)\setminus S(t+h)} \Big( {\pi \over 2} + \arctan u_x(t+h,\xi(t+h,y))\Big)
\,\phi_{1,0}(y)\,\bar u_x^2(y)\,dy$$
$$+\kappa_0 \int_{S(t) \setminus S(t+h)}  \Big( {\pi \over 2} +
\arctan \tilde u_x (t+h,\tilde \xi(t+h,\tilde y))\Big) \,\phi_{2,0}(\tilde y)\, \tilde u_x^2(0,\tilde y)\,\psi_0'(y)\,dy$$
$$- \int_{S(t) \setminus S(t+h)} E(t+h,y)\,\phi_{1,0}(y)\,\bar u_x^2(y)\,dy$$
$$=\kappa_0 \int_{S(t) \setminus S(t+h)} \Big( {\pi} + \arctan u_x(t+h,\xi(t+h,y))
+ \arctan \tilde u_x (t+h,\tilde \xi(t+h,\tilde y))\Big)
\,\phi_{1,0}(y)\,\bar u_x^2(y)\,dy$$
$$- \int_{S(t) \setminus S(t+h)} E(t+h,y)\,\phi_{1,0}(y)\,\bar u_x^2(y)\,dy\eqno(3.33)$$
by (3.10) for $(\phi_{1,0},\,\phi_{2,0},\, \psi_0)$. In view of (3.29)-(3.30),
on $S(t) \setminus S(t+h)$ we have that
$$E(t+h,y)={\pi} + \arctan u_x(t+h,\xi(t+h,y)) + \arctan \tilde u_x (t+h,\tilde \xi(t+h,\tilde y))$$
since at least one of the expressions $u_x(t+h,\xi(t+h,y))$ and
$\tilde u_x (t+h,\tilde \xi(t+h,\tilde y))$ is precisely $-\infty$ on this set. Thus the whole
expression (3.33) is identically zero. Therefore (3.32) yields
$$J(t+h)-J(t) = \int_{S(t)} \Big( E(t+h,y)-E(t,y)\Big) \,\phi_{1,0}(y)\,\bar u_x^2(y)\,dy$$
$$ + \kappa_0 \int_{S(t)} \Big( \arctan u_x(t+h,\xi(t+h,y))-\arctan u_x(t,\xi(t,y))\Big)
\,\Big(1-\phi_{1,0}(y)\Big)\,\bar u_x^2(y)\,dy$$
$$ +  \kappa_0 \int_{S(t)} \Big( \arctan \tilde u_x (t+h,\tilde \xi(t+h,\tilde y))-
\arctan \tilde u_x (t,\tilde \xi(t,\tilde y))\Big)
\,\Big(1-\phi_{2,0}(\tilde y)\Big)\, \tilde u_x^2(0,\tilde y)\,\psi_0'(y)\,dy$$
$$ + \kappa_0 \int_{S^c(t)} \Big( \arctan u_x(t+h,\xi(t+h,y))-\arctan u_x(t,\xi(t,y))\Big)
\,\bar u_x^2(y)\,dy$$
$$ +  \kappa_0 \int_{S^c(t)} \Big( \arctan \tilde u_x (t+h,\tilde \xi(t+h,\tilde y))-
\arctan \tilde u_x (t,\tilde \xi(t,\tilde y))\Big)
\, \tilde u_x^2(0,\tilde y)\,\psi_0'(y)\,dy.\eqno(3.34)$$
In view of (3.29), we have $S^c(t)=\{T(y) \le t\}\,\cup\{\tilde T(\tilde y)\le t\}$.
On the set $\{T(y) \le t\}$ we have
$\arctan u_x(t+h,\xi(t+h,y))=\arctan u_x(t,\xi(t,y))=-\infty$ so that in the fourth term
in (3.34) only the integral over $\{T(y)>t,\,\tilde T (\tilde y) \le t\}$ might have a nonzero
contribution. Thus the second and fourth terms in (3.34) combine to
$$\kappa_0 \int_{\{T(y) >t\}}\Big( \arctan u_x(t+h,\xi(t+h,y))-\arctan u_x(t,\xi(t,y))\Big)
\,\Big(1-\phi_{1,0}(y)\Big)\,\bar u_x^2(y)\,dy$$
$$+\kappa_0 \int_{\{\tilde T(\tilde y) \le t < T(y)\}} \Big( \arctan u_x(t+h,\xi(t+h,y))-\arctan u_x(t,\xi(t,y))\Big)
\,\phi_{1,0}(y)\,\bar u_x^2(y)\,dy.\eqno(3.35)$$
Similarly, the third and fifth terms combine to
$$\kappa_0 \int_{\{\tilde T(\tilde y) >t\}}\Big( \arctan \tilde u_x(t+h,\xi(t+h,\tilde y))-\arctan\tilde u_x(t,\xi(t,\tilde y))\Big)
\,\Big(1-\phi_{2,0}(\tilde y)\Big)\,\tilde u_x^2(0,\tilde  y)\,\psi_0'(y)\,dy$$
$$+\kappa_0 \int_{\{T(y) \le t <\tilde T(\tilde y)\}} \Big( \arctan \tilde u_x(t+h,\xi(t+h,\tilde y))-\arctan\tilde u_x(t,\xi(t,\tilde y))\Big)
\,\phi_{2,0}(\tilde y)\,\tilde u_x^2(0,\tilde  y)\,\psi_0'(y)\,dy.\eqno(3.36)$$
To transform suitably the first term in (3.34), let us denote by $E^1(t,y)$ the
first expression in the minimum (3.30), and by $E^2(t,y)$ the second. If $E(t,y)=E^2(t,y)$,
then
$$E(t+h,y)-E(t,y) \le \kappa_0\,\Big( \arctan u_x(t+h,\xi(t+h,y))-\arctan u_x(t,\xi(t,y)) $$
$$+ \arctan \tilde u_x(t+h,\tilde \xi(t+h,\tilde y)) - \arctan \tilde u_x(t,\tilde \xi(t,\tilde y))\Big),\eqno(3.37)$$
since $E(t+h,y) \le E^2(t+h,y)$. On the other hand, if $E(t,y)=E^1(t,y)$, then the triangle
inequality and the relation $E(t+h,y) \le E^1(t+h,y)$ ensure that
$$E(t+h,y)-E(t,y) \le \Big|\xi(t+h,y)-\xi(t,y)+\tilde \xi(t+h,\tilde y)- \tilde \xi(t,\tilde y)\Big| $$
$$+\Big| u(t+h,\xi(t+h,y))-\tilde u(t+h,\tilde \xi(t+h,\tilde y))-u(t,\xi(t,y))+
\tilde u(t,\xi(t,\tilde y))\Big| \eqno(3.38)$$
$$+\kappa_0\,\Big|\arctan u_x(t+h,\xi(t+h,y))-\arctan u_x(t,\xi(t,y))
+ \arctan \tilde u_x(t+h,\tilde \xi(t+h,\tilde y)) -
\arctan \tilde u_x(t,\tilde \xi(t,\tilde y))\Big| .$$
Letting $h \downarrow 0$ in (3.34), and taking into account (3.38) and the
considerations preceding  it, we deduce that
$$\limsup_{h \downarrow 0} {J(t+h)-J(t) \over h} \le \kappa_0 J_0(t)+
\int_{S(t)} \phi_{1,0}(y)\,\bar u_x^2(y)\, \Big| {d \over dt}\,\xi(t,y)-{d \over dt}\,\tilde \xi(t,\psi_0(y))\Big|\,dy$$
$$+ \int_{S(t)} \phi_{1,0}(y)\,\bar u_x^2(y)\, \Big| {d \over dt}\,u(t,\xi(t,y))-{d \over dt}\, \tilde u(\tilde \xi(t,\psi_0(y)))\Big|\,dy$$
$$+\kappa_0 \int_{S(t)} \phi_{1,0}(y)\,\bar u_x^2(y)\, \Big| {d \over dt}\,\arctan\,u_x(t,\xi(t,y))-
{d \over dt}\, \arctan\,\tilde u_x(\tilde \xi(t,\psi_0(y)))\Big|\,dy \eqno(3.39)$$
where
$$J_0(t)= \int_{[T(y)>t]} \Big(1-\phi_{1,0}(y)\Big)\,\Big[
{d \over dt}\,\arctan\,u_x(t,\xi(t,y))\Big]\,\bar  u_x^2(y)\,dy$$
$$+ \int_{[\tilde T(\tilde y)>t]} \Big(1-\phi_{2,0}(\tilde y)\Big)\,\Big[
{d \over dt}\,\arctan\,\tilde u_x(t,\xi(t,\tilde y))\Big]\,  \tilde u_x^2(0,\tilde y)\,d\tilde y$$
$$+ \int_{[\tilde T(\tilde y) \le t <T(y)]} \phi_{1,0}(y)\,\Big[
{d \over dt}\,\arctan\,u_x(t,\xi(t,y))\Big]\,\bar  u_x^2(y)\,dy$$
$$+ \int_{[T(y) \le t < \tilde T(\tilde y)]} \phi_{2,0}(\tilde y)\,\Big[
{d \over dt}\,\arctan\,\tilde u_x(t,\xi(t,\tilde y))\Big]\,  \tilde u_x^2(0,\tilde y)\,d\tilde y \le 0,$$
the last inequality being true by (2.4).

Before proceeding with the further analysis of (3.39), we establish a few
{\it a priori} bounds. From (2.22) we get
$$\Big| {d \over dt}\, \xi(t,y)-{d \over dt}\,\tilde \xi(t,\tilde y) \Big| =\Big| u(t,\xi(t,y))-
\tilde u(t,\tilde \xi(t,\tilde y))\Big|.\eqno(3.40)$$
Also, note that if $v=\arctan z(t)$ and $\dot z=-\,\displaystyle{z^2 \over 2}$, then
$$\dot v=-\,{z^2 \over 2+2z^2}=-\,{1 \over 2}\,\sin^2\,v.$$
Since $|\sin^2\alpha - \sin^2\beta| \le |\alpha - \beta|$ by the mean-value theorem as
$|(\sin^2 z)'|=2|\sin z\,\cos z|=|\sin\,(2z)| \le 1$, we infer three useful facts if we set
$z=u_x(t,\xi(t,x))$. First of all,
$${d \over dt}\, \arctan u_x(t,\xi(t,y)),\ {d \over dt}\, \arctan \tilde u_x(t,\tilde u(t,\tilde
\xi(t,\tilde y)) \le 0.\eqno(3.41)$$
Secondly,
$$\Big| {d \over dt}\,\arctan\,u(t,\xi(t,y))-{d \over dt}\,\arctan\,
\tilde u (t,\tilde \xi(t,\tilde y))\Big| \le {1 \over 2}\,
\Big| \arctan u_x(t,\xi(t,y)) - \arctan \tilde u_x(t,\tilde \xi(t,\tilde y))\Big|.\eqno(3.42)$$
Furthermore, if $\arctan z \le -\,\displaystyle{\pi \over 4}$, then $\sin(\arctan z) \in [-1,-\,
\displaystyle{1 \over \sqrt{2}}]$, so that
$${d \over dt} \arctan u_x(t,\xi(t,y))\le -\,{1 \over 4}\quad\hbox{if}\quad
\arctan u_x(t,\xi(t,y)) \le -\, {\pi \over 4}.\eqno(3.43)$$
On the other hand, using first (2.21) and then (2.18), we have
$$\Big| {d \over dt}\,u(t,\xi(t,y_0))-{d \over dt}\,\tilde u (t,\tilde \xi(t,\tilde y_0))\Big|
=\Big| \varphi(t,y_0)-\tilde \varphi(t,\tilde y_0)) \Big|$$
$$={1 \over 4}\,\Big| \int_{\{T(y)>t\}} \sgn(y_0-y)\,\bar u_x^2(y)\,dy
- \int_{\{ \tilde T(\tilde y) > t \}} \sgn(\tilde y_0 - \tilde y)\,\tilde u_x^2(0,\tilde y)\,d\tilde y\Big|$$
$$={1 \over 4}\,\Big| \int_{-\infty}^{y_0} \bar u_x^2(y)\chi_{[T(y)>t]}dy -\int_{-\infty}^{\tilde y_0}
\tilde u_x^2(0,\tilde y)\,\chi_{[\tilde T(\tilde y)>t]}d\tilde y$$
$$-\int_{y_0}^{\infty} \bar u_x^2(y)\chi_{[T(y)>t]}dy -\int_{\tilde y_0}^{\infty}
\tilde u_x^2(0,\tilde y)\,\chi_{[\tilde T(\tilde y)>t]}d\tilde y\Big|$$
$$={1 \over 4}\,\Big| \int_{-\infty}^{y_0} \bar u_x^2(y)\chi_{[T(y)>t]}dy -\int_{-\infty}^{y_0}
\tilde u_x^2(0,\psi_0(y))\,\chi_{[\tilde T(\tilde y)>t]}\,\psi_0'(y)\,dy$$
$$-\int_{y_0}^{\infty} \bar u_x^2(y)\chi_{[T(y)>t]}dy -\int_{y_0}^{\infty}
\tilde u_x^2(0,\psi_0(y))\,\chi_{[\tilde T(\tilde y)>t]}\,\psi_0'(y)\,dy\Big|$$
$$={1 \over 4}\,\Big| \int_{-\infty}^{y_0} \Big(1-\phi_{1,0}(y)+\phi_{1,0}(y)\Big)\bar u_x^2(y)\chi_{[T(y)>t]}dy$$
$$ -\int_{-\infty}^{y_0} \Big( 1-\phi_{2,0}(\psi_0(y))+\phi_{2,0}(\psi_0(y))\Big)
\tilde u_x^2(0,\psi_0(y))\,\chi_{[\tilde T(\tilde y)>t]}\,\psi_0'(y)\,dy$$
$$-\int_{y_0}^{\infty} \Big(1-\phi_{1,0}(y)+\phi_{1,0}(y)\Big)\bar u_x^2(y)\chi_{[T(y)>t]}dy $$
$$-\int_{y_0}^{\infty} \Big( 1-\phi_{2,0}(\psi_0(y))+\phi_{2,0}(\psi_0(y))\Big)
\tilde u_x^2(0,\psi_0(y))\,\chi_{[\tilde T(\tilde y)>t]}\,\psi_0'(y)\,dy\Big|$$
after performing in the next to the last step in two of the integrals the change
of variables $\tilde y=\psi_0(y)$. Since (3.10) for $(\psi_0,\phi_{1,0},\phi_{2,0})$
ensures that $\phi_{1,0}(y)\,\bar u_x^2(y)=\phi_{2,0}(\psi_0(y))\,\tilde u_x^2(0,\psi_0(y))\,\psi_0'(y)$ a.e. on $S(t)$, we deduce that
$$\Big| {d \over dt}\,u(t,\xi(t,y_0))-{d \over dt}\,\tilde u (t,\tilde \xi(t,\tilde y_0))\Big|$$
$$\le {1 \over 4}\,\int_{[T(y)>t]}\Big( 1 - \phi_{1,0}(y)\Big)\,\bar u_x^2(y)\,dy
+ {1 \over 4} \int_{[\tilde T(\tilde y)>t]} \Big( 1 - \phi_{2,0}(\tilde y)\Big) \tilde u_x^2(0,\tilde y)\,d\tilde y$$
$$+{1 \over 4}\,\int_{[\tilde T(\tilde y) \le t <T(y)]} \phi_{1,0}(y)\,\bar u_x^2(y)\,dy
+ {1 \over 4} \int_{[T(y) \le t <\tilde T(\tilde y)]}  \phi_{2,0}(\tilde y)\, \tilde u_x^2(0,\tilde y)\,d\tilde y.\eqno(3.44)$$
Let us now introduce the following sets
$$S_1=\{y \in \R:\ \arctan\,u_x(t,\xi(t,y)) \le -{\pi \over 4}\quad\hbox{and}\quad
\arctan\,\tilde u_x(t,\tilde \xi(t,\tilde y)) \le -{\pi \over 4}\},$$
$$S_2=\{y\in \R:\ \arctan\,u_x(t,\xi(t,y)) > -{\pi \over 4}\quad\hbox{and}\quad
\arctan\,\tilde u_x(t,\tilde \xi(t,\tilde y)) > -{\pi \over 4}\},$$
$$S_3=\{y\in \R:\ \arctan\,u_x(t,\xi(t,y)) > -{\pi \over 4} \ge
\arctan\,\tilde u_x(t,\tilde \xi(t,\tilde y))\},$$
$$S_4=\{y\in \R:\ \arctan\,\tilde u_x(t,\tilde \xi(t,\tilde y))> -{\pi \over 4} \ge \arctan\,u_x(t,\xi(t,y))\}.$$
The integral on the right-hand side of (3.44) over $S_1$ is, in view of (3.43),
bounded from above by $|J_0(t)|=-\,J_0(t)$. The integral over $S_2$ is, in view
of the formula for $J(t)$ preceding relation (2.29), bounded from above by
$\displaystyle{J(t) \over \pi\kappa_0}$. To evaluate the contribution over the integral
over $S_3$, notice that the same formula for $J(t)$ yields
$$\eqalign{J(t) &\ge \kappa_0\,\int_{S(t)\,\cap\,S_3}  \Big(\arctan u_x(t,\xi(t,y)) - \arctan \tilde u_x(t,\tilde \xi(t,\tilde y))\Big) \,\phi_{1,0}(y)\,\bar u_x^2(y)\,dy\cr
&+\kappa_0 \int_{S(t)\,\cap\,S_3} \Big( {\pi \over 2} + \arctan u_x(t,\xi(t,y))\Big)\,
\Big(1-\phi_{1,0}(y)\Big)\,\bar u_x^2(y)\,dy\cr
&+\kappa_0 \int_{S^c(t)\,\cap\,S_3}
\Big( {\pi \over 2} + \arctan u_x(t,\xi(t,y))\Big)\,\bar u_x^2(y)\,dy\cr
&+\kappa_0 \int_{S(t)\,\cap\,S_3} \Big( {\pi \over 2} + \arctan \tilde u_x(t,\tilde \xi(t,\,\tilde y))\Big)\,
\Big(1-\phi_{2,0}(\tilde y)\Big)\,\tilde u_x^2(0,\tilde y)\,\psi_0'(y)\,d y\cr
&+\kappa_0 \int_{S^c(t)\,\cap\,S_3}
\Big( {\pi \over 2} + \arctan \tilde u_x(t,\,\tilde \xi(t,\,\tilde y))\Big)\,
\tilde u_x^2(0,\tilde y)\,\psi_0'(y)\,dy.\cr}$$
Using (3.10), the sum of the first term and the fourth term is larger than
$$\kappa_0\int_{S(t)\,\cap\,S_3}  \Big( {\pi \over 2} + \arctan u_x(t,\xi(t,y))\Big)\,\tilde u_x^2(0,\tilde y)\,\psi_0'(y)\,dy $$
$$\ge {\kappa_0\pi \over 4}\int_{S(t)\,\cap\,S_3}
\tilde u_x^2(0,\tilde y)\,\psi_0'(y)\,d y={\kappa_0\pi \over 4}\int_{S(t)\,\cap\,S_3}
\tilde u_x^2(0,\tilde y)\,d \tilde y$$
$$\ge {\kappa_0\pi \over 4}\int_{S(t)\,\cap\,S_3} \Big(1-\phi_{2,0}(\tilde y)\Big)\,
\tilde u_x^2(0,\tilde y)\,d \tilde y=
{\kappa_0\pi \over 4}\int_{[\tilde T(\tilde y)>t]\,\cap\,S_3} \Big(1-\phi_{2,0}(\tilde y)\Big)\,
\tilde u_x^2(0,\tilde y)\,d \tilde y,$$
with the last equality enforced by $S_3 \subset [T(y)>t]$.
The second term is bounded from below by
$${\kappa_0\pi \over 4} \int_{S(t) \,\cap\,S_3}\Big(1-\phi_{1,0}(y)\Big)\,\bar u_x^2(y)\,dy
={\kappa_0 \pi \over 4} \int_{[T(y)>t]\,\cap\,S_3}\Big(1-\phi_{1,0}(y)\Big)\,\bar u_x^2(y)\,dy$$
$$-{\kappa_0\pi \over 4} \int_{[T(y)>t \ge \tilde T(\tilde y)]\,\cap\,S_3}
\Big(1-\phi_{1,0}(y)\Big)\,\bar u_x^2(y)\,dy={\kappa_0 \pi \over 4} \int_{[T(y)>t]\,\cap\,S_3}\Big(1-\phi_{1,0}(y)\Big)\,\bar u_x^2(y)\,dy$$
$$-{\kappa_0\pi \over 4} \int_{[T(y)>t \ge \tilde T(\tilde y)]\,\cap\,S_3} \bar u_x^2(y)\,dy + {\kappa_0\pi \over 4} \int_{[T(y)>t \ge \tilde T(\tilde y)]\,\cap\,S_3} \phi_{1,0}( y)\,\bar u_x^2(y)\,dy$$
if we recall (3.10).
The third term is bounded from below by
$${\kappa_0\pi \over 4} \int_{S^c(t) \,\cap\,S_3}  \bar u_x^2(y)\,dy
\ge {\kappa_0\pi \over 4} \int_{[T(y)>t \ge \tilde T(\tilde y)]\,\cap\,S_3} \bar u_x^2(y)\,dy.$$
Summing up, we get
$$J(t) \ge {\kappa\pi \over 4}\Big\{\int_{[T(y)>t]\,\cap S_3}\Big( 1 - \phi_{1,0}(y)\Big)\,\bar u_x^2(y)\,dy
+  \int_{[\tilde T(\tilde y)>t]\,\cap S_3} \Big( 1 - \phi_{2,0}(\tilde y)\Big) \tilde u_x^2(0,\tilde y)\,d\tilde y$$
$$+\,\int_{[\tilde T(\tilde y) \le t <T(y)]\,\cap S_3} \phi_{1,0}(y)\,\bar u_x^2(y)\,dy
+ \int_{[T(y) \le t <\tilde T(\tilde y)]\,\cap S_3}  \phi_{2,0}(\tilde y)\, \tilde u_x^2(0,\tilde y)\,d\tilde y\Big\},$$
since the last term on the right is zero as $S_3 \subset [T(y)>t]$. A similar relation holds with
$S_4$ instead of $S_3$. Consequently, putting together all this information about
the various inequalities valid on the disjoint sets $S_1,\,S_2,\,S_3,$ and $S_4$, we
conclude by (3.44) that
$$\Big| {d \over dt}\,u(t,\xi(t,y_0))-{d \over dt}\,\tilde u (t,\tilde \xi(t,\tilde y_0))\Big|
\le -J_0(t)+{3J(t) \over \kappa_0 \pi}.\eqno(3.45)$$

To obtain  now a suitable estimate on
$$\limsup_{h \downarrow 0}{J(t+h)-J(t) \over h} =
\limsup_{h \downarrow 0}\int_{S(t)} {E(t+h,y)-E(t,y) \over h}\,\phi_{1,0}(y)\,\bar u_x^2(y)\,dy + \kappa_0 \,J_0(t)\eqno(3.46)$$
we distinguish two cases. If $E(t,y)$ is the
second component $E^2(t,y)$ of the minimum in (3.30), then by (3.37)
and (3.41) we can estimate
the contribution of the first term in (3.46) by zero from above. On the other hand,
if the minimum is $E^1(t,y)$, then the first integral term in (3.46) is not larger (pointwise)
than the nonnegative expression
$$\Big(E(t,y) + {3J(t) \over \kappa_0 \pi} -J_0(t)\Big)\,\phi_{1,0}(y)\,\bar u_x^2(y)$$
in view of the estimates (3.40), (3.42), and (3.45). We conclude that
$$\limsup_{h \downarrow 0}{J(t+h)-J(t) \over h} \le J(t) + \Big({3J(t) \over  \kappa_0\pi} -J_0(t)
\Big)\,\Vert \bar u_x^2 \Vert_{L^1(\R)} +\kappa_0 \,J_0(t).$$
Since $J_0(t) \le 0$, choosing the constant $\kappa_0\doteq \|\bar u^2_x\big\|_{\L^1(\R)}$
we now have
$${d\over dt} J^{(\psi^t,\phi_1^t,\phi_2^t)}
\big(u(t),\,v(t)\big)\leq 2\,J^{(\psi^t,\phi_1^t,\phi_2^t)}
\big(u(t),\,v(t)\big)\,.$$
Optimizing over all triples $(\psi^0,\phi_1^0,\phi_2^0)$
we conclude
$$J\big( u(t),\, v(t)\big) \leq
J\big(u(0),\, v(0)\big)\,e^{2t},\quad t \ge 0.\eqno(3.47)$$

\v
Summing up the considerations made above, we proved the
following result.

\v
{\bf Theorem 2.} {\it The trajectories $t \mapsto u(t)$ of $(1.1)$
constructed in Theorem 1 are locally Lipschitz continuous as maps from $[0,\infty)$ into the metric space ${\cal X}$ equipped with the
distance functional $J$. Moreover, the distance between two
trajectories is also locally Lipschitz continuous as a map from
$[0,\infty)$ into ${\cal X}$.}

\vsk
\n{\medbf 4 - Concluding remarks}
\v
The following example shows that, in some sense, our distance
functional
$J$ in (3.11) is ``sharp''.  Indeed, the convergence
of the initial data in  $\L^\infty(\R)\cap \L^1(\R)$ together with the
weak convergence of the derivatives $\bar u_x$ and
$\bar u_x^2$ in $\L^2(\R)$ does not
guarantee the convergence of the corresponding solutions
at later times $t>0$.
\v
\n{\bf Example 1.}
Consider the functions $f,g:[0,1]\mapsto [0,1]$
defined as
$$f(x)\doteq\cases{ 1-2x\qquad &if\quad $x\in [0,~1/2]$,\cr
0&if\quad $x\in [1/2\,,~1]$,\cr}
\qquad\qquad g(x)\doteq\cases{ 1-3x\qquad &if\quad $x\in [0,~1/6]$,\cr
1/2\qquad &if\quad $x\in [1/6\,,~1/2]$,\cr
1-x\qquad &if\quad $x\in [1/2\,,~1]$.\cr}$$
Observe that
$$\int_0^1 f'(x)\,dx=\int_0^1 g'(x)\,dx=-1\,,\qquad\qquad
\int_0^1 \big[f'(x)\big]^2\,dx=\int_0^1 \big[g'(x)\big]^2\,dx=2\,.$$
Next, consider the function
$$h(x)\doteq\cases{ 1-|x|\qquad &if\quad $|x|\leq 1$,\cr
0&if\quad $|x|\geq 1$,\cr}$$
and define the sequences of initial values
$$\bar u_n(x)=\cases{h(x)\qquad &if\quad $x\notin [0,1]$,\cr
h(i/n)+{1\over n}
f(nx-i+1)\qquad &if\quad $x\in \left[ {i-1\over n}\,,~{i\over n}\right]
\qquad\quad i=1,\ldots,n$,
\cr}$$
$$\bar v_n(x)=\cases{h(x)\qquad &if\quad $x\notin [0,1]$,\cr
h(i/n)+{1\over n}
g(nx-i+1)\qquad &if\quad $x\in \left[ {i-1\over n}\,,~{i\over n}\right]
\qquad\quad i=1,\ldots,n$.
\cr}$$
Letting $n\to\infty$ we now have the strong convergence
$\big\|\bar u_n-\bar v_n\big\|_{\L^\infty(\R)}\to 0$.
Moreover, by construction it is easy to see that at each point
$x \in [0,1]$,
$$\lim_{n \to \infty} \int_0^x \Big( (\bar u_n)_x(y)-(\bar v_n)_x(y)\Big)\,dy=
\lim_{n \to \infty} \int_0^x \Big( (\bar u_n)^2_x(y)-(\bar v_n)^2_x(y)\Big)\,dy
=0$$
so that in $\L^2[0,1]$ one has the weak convergence
$$(\bar u_n)_x-(\bar v_n)_x\wto 0\,,\qquad\qquad (\bar u_n)^2_x-(\bar
v_n)^2_x
\wto 0\,,$$
since both sequences are bounded in $\L^2[0,1]$ and the previous
observation identifies the zero function as the only possible weak limit. However
$$u(t)\doteq \lim_{n\to\infty} u_n(t)~\not=~ \lim_{n\to\infty}v_n(t)
\doteq v(t)$$
for every $t \in (2/3,1)$, where
$T=2/3$ is the time at which the gradients of the functions
$v_n$ blow up. The last assertion follows at once from (2.27).$\hfill \diamondsuit$

\vs
We also would like to highlight the importance of
requiring that the transport map
$\psi$ in (3.10) be monotone nondecreasing.
If in (3.11) we were to take the minimization over all maps
$\psi$, not necessarily monotone, we would obtain the classical
Kantorovich-Rubinstein
distance between measures, which generates the weak topology
on the space of bounded, positive measures [V].
By restricting ourselves to monotone nondecreasing
maps $\psi$, the corresponding distance functional
generates a much stronger topology.
\v
\n{\bf Example 2.} Consider the sequence of Lipschitz functions
$$u^m(x)\doteq\cases{ 0\qquad &if\quad $x\notin [0,1]$\cr
x-{(i-1)/m}\qquad &if\quad $(i-1)/m\leq x\leq (2i-1)/2m$\cr
{i/m}-x\qquad &if\quad $(2i-1)/2m\leq x\leq i/m$\cr}
\qquad i=1,\ldots,m\,.$$
In this case, $u_x=\pm 1$ and $\arctan u_x=\pm \pi/4$.
The corresponding measures $\mu^{u^m}$ defined at (3.9)
converge weakly to the measure
$\mu$ on $\R^2\times [-\pi/2\,,~\pi/2]$ defined as
$$\mu(A)\doteq {1\over 2}\,\meas\Big\{ x\in [0,1]\,;~~
(x,\,0,\,\pi/4)\in A\Big\}
+{1\over 2}\,\meas\Big\{ x\in [0,1]\,;~~(x,\,0,\,-\pi/4)\in A\Big\}\,.$$
In particular, these measures form a
Cauchy sequence in the Kantorovich metric. However, these same
functions $u^m$ do not form a Cauchy sequence w.r.t.~the distance
$J$. Indeed, let $m<n$.  Consider the open intervals
$$I^{m+}_i=\,\left]{i-1\over m}\,,{2i-1\over 2m}\right[\,,
\qquad\qquad
I^{m-}_i=\,\left]{2i-1\over 2m}\,,{i\over m}\right[\,,
$$
where $u_x^m$ takes the values $+1$ and $-1$, respectively.
Define the intervals $I^{n+}_j$, $I^{n-}_j$ similarly.
Now consider any transportation plan
$(\psi,\phi_1,\phi_2)$, with $\psi$ non-decreasing.
For each $i=1,\ldots,m$,
call $\nu_i$ the number of distinct intervals $I^{n+}_j$
which intersect the image $\psi(I^{m+}_i)$.
Since $\psi$ is monotone, if $\nu_i\geq 2$, this implies that
the image $\psi(I^{m+}_i)$ entirely covers $\nu_i-1$
distinct intervals $I^{n-}_j$.  Because $u^m_x=1$ on $I^{m+}_i$
and $u^n_x=-1$ on each $I^{n-}_j$, this accounts for a cost
$\geq (\nu_i-1)/2n $.
Next, if $\nu_1+\ldots +\nu_m =n^* < n$, there must be $n-n^*$ intervals
$I^{n+}_{j(1)}\,,\ldots, \, I^{n^+}_{j(n^*-n)}$
which do not intersect any of the sets $\psi(I^{m+}_i)$, for $i=1,\ldots,m$.
These intervals must be contained in the image of some $I^{m-}_i$,
or in the image of the set $\psi\big(\R\setminus [0,1]\big)$,
where $u^m\equiv 0$.
This accounts for a cost $\geq (n-n^*)/4n $.

The above argument shows that, for any $m<n$, the cost of any
transportation plan is bounded below by
$$J^{(\psi,\phi_1,\phi_2)} (u^m, u^n)\geq
{1\over 4n}\cdot\max\left\{ \sum_{i=1}^m (\nu_i-1)\,,~n-\sum_{i=1}^m \nu_i
\right\}
\geq {1\over 4n}\cdot {n-m\over 2}\,.$$
For any fixed $m$, the right hand side  approaches $1/8$ as
$n\to\infty$. Therefore, the above is not a Cauchy sequence,
in our transportation metric.
$\hfill \diamondsuit$

\v
{\it Acknowledgements.} A. Bressan was supported by the Italian M.I.U.R.,
within the research project 2002017219  "Equazioni iperboliche e
paraboliche non lineari" and A. Constantin was supported by the
Science Foundation Ireland within the research project 04/BR/M0042 "Nonlinear Waves".

\vsk

\c{\medbf References}
\v
\item{[AA]}
G.~Alberti and L.~Ambrosio,
A geometrical approach to monotone functions in $\R^n$,
{\it Math. Z.} {\bf 230} (1999), 259-316.
\v
\item{[BSS]}
R.~Beals, D.~Sattinger, and J.~Szmigielski,
Inverse scattering solutions of the Hunter-Saxton equation,
{\it Applicable Analysis} {\bf 78} (2001), 255-269.
\v
\item{[B]}
A.~Bressan,
Unique solutions for a class of discontinuous differential equations,
{\it Proc. Amer. Math. Soc.} {\bf 104} (1988), 772-778.
\v
\item{[BC]}
A.~Bressan and G.~Colombo,
Existence and continuous dependence for discontinuous O.D.E.'s,
{\it Boll. Un. Mat. Ital.} {\bf 4-B} (1990), 295-311.
\v
\item{[BZZ]}
A.~Bressan, P.~Zhang, and Y. Zheng,
On asymptotic variational wave equations,
{\it Arch. Rat. Mech. Anal.}, to appear.
\v
\item{[BGH]}
G.~Buttazzo, M.~Giaquinta, and S.~Hildebrandt,
{\it One-dimensional Variational Problems},
Clarendon Press, Oxford, 1998.
\v
\item{[CH]}
R.~Camassa and D.~D.~Holm,
An integrable shallow water equation with peaked solitons,
{\it Phys. Rev. Lett.}
{\bf 71} (1993), 1661-1664.
\v
\item{[C]}
A.~Constantin,
Existence of permanent and breaking waves for a shallow water
equation: a geometric approach, {\it Ann. Inst. Fourier (Grenoble)}
{\bf 50} (2000), 321-362.
\v
\item{[CE]}
A.~Constantin and J.~Escher,
Wave breaking for nonlinear nonlocal shallow water equations,
{\it Acta Mathematics} {\bf 181} (1998), 229-243.
\v
\item{[CK]}
A.~Constantin and B.~Kolev,
Geodesic flow on the diffeomorphism group of the circle,
{\it Comment. Math. Helv.} {\bf 78} (2003),
787-804.
\v
\item{[DP]}
H.~H.~Dai and M.~Pavolv,
Transformations for the Camassa-Holm equation, its high-frequency
limit and the Sinh-Gordon equation,  {\it J. Phys. Soc. Japan}
{\bf 67} (1998), 3655-3657.
\v
\item{[EG]}
L.~C.~Evans and R.~F.~ Gariepy, {\it
Measure Theory and Fine Properties of Functions}, CRC Press, Boca Raton, FL, 1992.
\v
\item{[HS]}
J.~K.~Hunter and R.~A.~Saxton,
Dynamics of director fields, {\it SIAM J. Appl. Math.} {\bf 51} (1991),
1498-1521.
\v
\item{[HZ1]}
J.~K.~Hunter and Y.~Zheng,
On a completely integrable nonlinear hyperbolic variational equation,
{\it Physica D} {\bf 79} (1994), 361-386.
\v
\item{[HZ2]}
J.~K.~Hunter and Y.~Zheng,
On a nonlinear hyperbolic variational equation I.
{\it Arch. Rat. Mech. Anal.} {\bf 129} (1995), 305-353.
\v
\item{[HZ3]}
J.~K.~Hunter and Y.~Zheng,
On a nonlinear hyperbolic variational equation II.
{\it Arch. Rat. Mech. Anal.} {\bf 129} (1995), 355-383.
\v
\item{[J]}
R.~S.~Johnson,
Camassa-Holm, Korteweg-de Vries and related models for water waves, {\it
J. Fluid Mech.} {\bf 455} (2002), 63-82.
\v
\item{[KM]}
B.~Khesin and G.~Misiolek,
Euler equations on homogeneous spaces and Virasoro orbits,
{\it Adv. Math.}
{\bf 176} (2003), 116-144.
\v
\item{[K]}
S.~Kouranbaeva,
The Camassa-Holm equation as a geodesic flow on the diffeomorphism group,
{\it J. Math. Phys.} {\bf 40} (1999), 857-868.
\v
\item{[N]}
J.~Natanson, {\it Theory of Functions of a Real Variable}
F.~Ungar Publ., New York, 1964.
\v
\item{[V]}
C.~Villani, {\it Topics in Optimal Transportation},
Amer. Math. Soc., Providence 2003.
\v
\item{[Z]}
E.~Zeidler, {\it Nonlinear Functional Analysis and its Applications},
Springer-Verlag, New York, 1990.
\v
\item{[ZZ1]}
P.~ Zhang and Y.~Zheng,
On oscillations of an asymptotic equation of a nonlinear variational
wave equation,
{\it Asympt. Anal.} {\bf 18} (1998), 307-327.
\v
\item{[ZZ2]}
P.~ Zhang and Y.~Zheng,
On the existence and uniqueness of solutions to an asymptotic equation
of a variational wave equation, {\it Acta Math. Sinica}
{\bf 15} (1999), 115-130.
\v
\item{[ZZ3]}
P.~ Zhang and Y.~Zheng,
Existence and uniqueness of solutions of an asymptotic equation arising
from a variational wave equation with general data,
{\it Arch. Rat. Mech. Anal.} {\bf 155} (2000), 49-83.

\bye